\documentclass[11pt,a4paper,reqno]{amsart}

\usepackage{mathrsfs,amssymb}

\newtheorem{theorem}{Theorem}[section]
\newtheorem{lemma}[theorem]{Lemma}
\newtheorem{proposition}[theorem]{Proposition}
\newtheorem{corollary}[theorem]{Corollary}

\newtheorem{question}[theorem]{Question}

\theoremstyle{definition}

\theoremstyle{plain}

\newcommand{\BN}{\mathbb{N}}
\newcommand{\OpPref}{\operatorname{Pref}}
\newcommand{\OpFact}{\operatorname{Fact}}
\newcommand{\OpInt}{\operatorname{Int}}
\newcommand{\OpSuff}{\operatorname{Suff}}

\begin{document}
    \title[Automatic {R}ees matrix semigroups]{Automatic {R}ees matrix semigroups \\ over categories}
    \subjclass[2000]{20M05, 20M50, 18B40}
    \keywords{Rees matrix semigroups, categories, automatic structures}

\maketitle

\begin{center}
    Mark Kambites \\

    \medskip

    Fachbereich Mathematik / Informatik,  Universit\"at Kassel \\
    34109 Kassel, Germany \\

    \medskip

    \texttt{kambites@theory.informatik.uni-kassel.de} \\
\end{center}

\bigskip

\begin{abstract} We consider the preservation of the properties of
automaticity and prefix-automaticity in Rees matrix semigroups over
semigroupoids and small categories.  Some of our results are new or
improve upon existing results in the single-object case of Rees
matrix semigroups over semigroups.
\end{abstract}

\section{Introduction}\label{sec:intro}

In recent years, one of the most productive areas of combinatorial
group theory has been the theory of \textit{automatic groups}.
The expression of a finitely generated group using an \textit{automatic
structure} provides a basis for efficiently performing computations in
the group which may be hard or impossible given only a presentation.
Groups which admit automatic structures also display a number of
remarkable algebraic and geometric properties. \cite{Epstein92}

The usual language-theoretic definition of automaticity in groups lends
itself naturally to application in wider contexts. It is observed in
\cite{Epstein92} that the notion of automaticity extend naturally from
groups to \textit{groupoids}. More recently, Hudson \cite{Hudson96} has
introduced a definition of automaticity in
semigroups and monoids and a theory of
automatic semigroups has begun to emerge
\cite{Campbell99,Campbell01,Campbell00a,Descalco01,Duncan99,Hoffmann00a,
Hoffmann01,Hoffmann02,Hoffmann00b,Hoffmann01a}. In \cite{KambitesAutoCat}
we developed a common framework for these generalisations, in the form of
a theory of automaticity for \textit{small categories} and
\textit{semigroupoids}. We explored what could be learnt about automatic
small categories and semigroupoids, by applying the existing theory of
automatic semigroups. Semigroupoids and small categories play a vital
role in the structural theory of semigroups, so it also seems natural to ask,
conversely, whether automatic categories can tell us anything about
automatic semigroups.

A key recurring theme in structural semigroup theory is that of the
\textit{Rees matrix construction} (see, for example, \cite{Howie95}).
A number of interesting results have been proved concerning the relationship
between automaticity properties and Rees matrix constructions over groups
\cite{Campbell00a} and more generally over semigroups \cite{Descalco01,Silva00b}.
At the same time, Lawson \cite{Lawson00} has applied to the study of
abundant semigroups a form of Rees matrix construction over
semigroupoids. In \cite{Kambites05} we explored combinatorial aspects of
this construction showing that, under certain assumptions, combinatorial
properties such as finite generation and finite presentability are
preserved. In this paper, we consider the relationship
between automaticity and these constructions. Some of our results are new,
or improve upon existing results, even when specialised to the single-object
case of Rees matrix semigroups over semigroups.

In addition to this introduction, this paper comprises four sections. In
Section~\ref{sec:background}, we briefly recall the notions of generators
and relations for partial algebras which we studied in \cite{Kambites05},
and the definitions of regular path languages and of automatic and
prefix-automaticity introduced in \cite{KambitesAutoCat}. We also recall
some key results from \cite{KambitesAutoCat} which will be applied in
this paper.
In Section~\ref{sec:technical}, we prove some technical results concerning
automaticity in small categories and semigroupoids, which will be needed
in the sections that follow.

Section~\ref{sec:autorees} considers the relationship between
automaticity and Rees matrix constructions with zero over
semigroupoids. We show that a finitely generated Rees matrix semigroup
over an automatic semigroupoid is always automatic, and provide some
sufficient conditions for the underlying semigroupoid of an automatic
Rees matrix semigroup to be automatic. In Section~\ref{sec:pautorees}
we prove some related results for prefix-automaticity.
Section~\ref{sec:reeswithoutzero} extends these results to the case of
Rees matrix constructions without zero. Finally, Section~\ref{sec:remarks}
contains some remarks and open questions.

\section{Background}\label{sec:background}

In this section, we briefly recall a number of definitions and results
from \cite{Kambites05} and \cite{KambitesAutoCat}. For a more detailed
introduction, the reader should consult those papers.

By a \textit{(directed) graph} $X$ we mean a collection $X^0$ of \textit{vertices}
together with a collection $X^1$ of \textit{edges} and two functions
$\alpha, \omega : X^1 \to X^0$
 which determine respectively the
\textit{source}
and \textit{target} of each edge.

A \textit{path} $\pi$ in $X$ is a finite sequence $e_1 e_2 \dots e_n$ of (not
necessarily distinct) edges in $X$ such that $e_i \omega = e_{i+1} \alpha$
for $1 \leq i < n$. We define $\pi \alpha = e_1 \alpha$,
$\pi \omega = e_n \omega$,
and call these the \textit{source} and
\textit{target} respectively of the path $\pi$.
The \textit{length}
of the path $\pi$ is the number $n$ of edges; we denote it by $|\pi|$.
There is also a distinct \textit{empty path} of
length $0$ at each vertex, with source and target that vertex.
We identify each path of length $1$ with its single edge, and each vertex
with the empty path at that vertex. Thus, $X^0$ and $X^1$ are the sets of
paths in $X$ of length $0$ and of length $1$ respectively. We extend this
notation by writing $X^n$ for the set of paths of length $n$ in $X$. 
We also write $X^{\geq n}$ for the set of paths in $X^*$ of length $n$ or
more; $X^{> n}$, $X^{\leq n}$ and $X^{< n}$ are defined analogously in
the obvious way, including empty paths where appropriate.

Let $X$ and $Y$ be graphs. The \textit{direct product} $X \times Y$ of $X$ and
$Y$ is the (finite)
graph with
vertex set $X^0 \times Y^0$, and edge set $X^1 \times Y^1$
$(f,g) \alpha = (f \alpha, g \alpha)$ and $(f,g) \omega = (f \omega, g \omega)$
for all $f \in X^1, g \in Y^1$. A \textit{graph morphism} $\sigma : X \to Y$ consists of
functions $\sigma^0 : X^0 \to Y^0$ and $\sigma^1 : X^1 \to Y^1$ such
that $x \sigma^1 \alpha = x \omega \alpha^0$ and
$x \sigma^1 \omega = x \omega \sigma^0$ for all edges $x \in X^1$.
If, in addition, $X^0 = Y^0$ and $\sigma^0 : X^0 \to X^0$ is the identity
function, then we call
$\sigma$ an \textit{edge-morphism}, and, for notational convenience, identify
$\sigma$ with $\sigma^1$.

A \textit{semigroupoid} $S$ is a small graph $S$,
together with a partial multiplication
on the edges of $S$ such that, for any two edges $e, f \in S^1$
\begin{itemize}
\item[(i)] the product $ef$ is defined if and only if $e \omega = f \alpha$;
\item[(ii)] if $ef$ is defined, then $(ef)\alpha = e \alpha$ and
            $(ef)\omega = f \omega$; and
\item[(iii)] if $ef$ and $fg$ are defined then $e(fg) = (ef)g$.
\end{itemize}
The vertices and edges of a semigroupoid are called respectively
\textit{objects} and \textit{arrows}.
Where no ambiguity arises, we abuse the notation slightly by writing $S$
to mean the set $S^1$ of arrows in $S$.

If $s_0 \in S^0$ is such that there exists an $s \in S^1$ with
$s \alpha = s \omega = s_0$,
we say that the \textit{local semigroup} of $S$ at
the object $s_0$ is the semigroup with elements
$\lbrace s \in S \mid s \alpha = s \omega = s_0 \rbrace$
and multiplication defined by restricting that in $S$.
 A semigroupoid
is \textit{strongly connected} if
for every pair $s_0, t_0 \in S^0$ of objects, there is an arrow (or
equivalently, a non-empty path) in $S$ with source $s_0$ and target $t_0$.
An object $x \in S^0$ of a semigroupoid $S$ is 
\textit{isolated} if is not the source or target of any arrow; if a
semigroupoid has no isolated objects then it is \textit{isolation-free}.

The \textit{free semigroupoid $X^+$} on a small graph $X$ is 
the semigroupoid whose objects are the vertices of $X$, and whose arrows
are the \textit{non-empty} paths in $X$, with $\alpha$, $\omega$ and the
partial multiplication defined in the obvious way. The \textit{free
category $X^*$} is the category obtained by adjoining a local identity
(which can be thought of as the empty path) at each object of $X^+$. A
\textit{path language} or just a \textit{language} over $X$ is a (possibly
empty) collection of (possibly empty) paths in $X$, that is, a subset of
the free category $X^*$.

Let $x = x_1 \dots x_n$ be a path in a free semigroupoid $X^+$ where each
$x_i \in X^1$. A \textit{prefix} of $x$ is a (non-empty) path of the form
$x_1 \dots x_j$ for some $1 \leq j \leq n$. A \textit{suffix} of $x$ is a
(non-empty) path of the form $x_j \dots x_n$ for some
 $1 \leq j \leq n$. A \textit{factor} of $x$ is a (non-empty) path of the
form $x_j \dots x_k$ for some $1 \leq j \leq k \leq n$. An   
\textit{internal factor} of $x$ is a (non-empty) path of the form $x_j
\dots x_k$ for some $1 < j \leq k < n$, that is, a factor of $x_2 \dots
x_{n-1}$.

Given a path language $L \subseteq X^+$, we write $\OpPref(L)$
[respectively $\OpSuff(L)$, $\OpFact(L)$, $\OpInt(L)$] for the set of
\textbf{non-empty} prefixes [respectively suffixes, factors, internal factors]
of paths in $L$. For $n \in \BN$ we write $\OpPref_n(L)$ to denote
$\OpPref(L) \cap X^n$, and similarly for $\OpSuff_n(L)$, $\OpFact_n(L)$  
and $\OpInt_n(L)$; note that all four are empty when $n = 0$.

A \textit{(path) automaton} $M$ over a
small graph $X$ is a small graph $M$, together with a graph morphism     
$\sigma : M \to X$, a set of distinguished \textit{start vertices} of $M$
and a set of distinguished \textit{terminal vertices} of $M$.
We think of the automaton as the graph $M$ with each edge $e \in M^1$
labelled by $e \sigma^1 \in X^1$, and each vertex $v \in M^0$ labelled by
$v \sigma^0 \in X^0$. The vertices and edges are called \textit{states}
and \textit{transitions} respectively. The unique \textit{language
accepted} or \textit{language recognised} by $M$ is the set of paths in
$X$ which label paths from a start vertex to a terminal vertex in $M$.

The automaton is called a \textit{complete, deterministic}
automaton
if (i) there is exactly one start vertex in the pre-image of each object
           in $X$ and
(ii) for every state $m \in M^0$ and every edge
$e \in X^1$ with $e \alpha = m \sigma^0$ there is a unique edge $f \in M^1$
with $f \sigma^1 = e$ and $f \alpha = m$.
The automaton is called \textit{finite}
if the graphs $M$ and $X$ are finite; we shall be concerned exclusively with
finite automata.

A path language which is accepted by some finite path automaton is called \textit{regular}.
We recall from \cite[Section~3]{KambitesAutoCat} that a language
$L \subseteq X^+ \subseteq (X^1)^+$ is regular in this sense if and only
if it is regular in the usual sense as a language over the alphabet $X^1$.
We recall also that the set of regular path languages over $X$ contains
$X^+$, $X^*$ and all finite path languages, and is closed under
concatenation, finite intersection, finite union, complement, set
difference, generation of subcategories and subsemigroupoids, and
prefix-closure. We will use these properties without further comment.

Given a graph $X$, we denote by $X^\$$ the graph $X$ with an extra edge $\$_s$
adjoined for every vertex $s \in X^0$, with source and target $s$. We define a
function $\delta_X : (X^+ \times X^+) \to (X^\$ \times X^\$)^+$ by
$$(a, b) \delta_X =
\begin{cases}
(a_1, b_1) \dots (a_n, b_n) & \text{ if } m = n \\
(a_1, b_1) \dots (a_n, b_n) (a_{n+1}, \$_{b_n \omega}) \
  \dots (a_m, \$_{b_n \omega}) & \text{ if } m > n \\
(a_1, b_1) \dots (a_m, b_m) (\$_{a_m \omega}, b_{n+1}) \
  \dots (\$_{a_m \omega}, b_n) & \text{ if } n > m \\
\end{cases}
$$
where $a = a_1 \dots a_m$ and $b = b_1 \dots b_n$.
We observe that the function $\delta_X$ is injective, and in particular
that it distributes over intersection, that is, that 
$$(R_1 \cap R_2) \delta_X = R_1 \delta_X \cap R_2 \delta_X$$
for all binary relations $R_1, R_2 \subseteq X^+ \times X^+$.

A \textit{(binary) synchronous path automaton} over $X$ is a
finite path automaton over the graph $X^\$ \times X^\$ $, with the property that
the language accepted is contained within the image
$(X^+ \times X^+) \delta_X$ of $\delta_X$.
For convenience, we shall say that a synchronous path automaton accepts a pair
$(a,b) \in X^+ \times X^+$
if it accepts $(a,b) \delta_X$. A binary relation    
$R \subseteq X^+ \times X^+$ is called \textit{synchronously regular} if there
exists a synchronous path automaton accepting exactly the language $R \delta_X$,
that is, if $R \delta_X$ is regular.
\begin{lemma}\cite[Lemma~3.7]{KambitesAutoCat}\label{lem:syncrat_prop}
Let $X$ be a finite graph. Then:
\begin{itemize}
\item[(i)] If $R \subseteq X^+ \times X^+$ is synchronously regular, then
the binary relation
$$R^{-1} = \lbrace (v,u) \mid (u, v) \in R \rbrace$$
is synchronously regular.
\item[(ii)] If $R \subseteq X^+ \times X^+$ is synchronously regular, then
the projections $R \pi_1$
and $R \pi_2$ of $R$ onto its first and second coordinates are regular.
\item[(iii)] If $K, L \subseteq X^+$ are regular languages of non-empty
paths, then the binary relation $K \times L$ is synchronously regular.
\item[(iv)] Synchronously regular binary relations over $X$ are closed under
 finite union, finite intersection and relational composition.
\item[(v)] If $K \subseteq X^+$ is regular then the \textit{diagonal binary
relation}
$$\lbrace (w, w) \mid w \in K \rbrace$$
is synchronously regular.
\end{itemize}
\end{lemma}

Let $X$ and $Y$ be finite graphs, $A \subseteq X^*$ and $\phi : A \to Y^*$
be a function.
We say that $\phi$ is \textit{strongly regularity preserving} if for  
every regular language $L \subseteq X^+$, we have that
$(L \cap A) \phi \subseteq Y^+$ is a regular language.

Let $X$ and $Y$ be finite graphs, and for $i = 1,2$ suppose we have
subsets $A_i \subseteq X^+$ and functions $\phi_i : A_i \to Y^+$.
We say that $\phi_1$ and $\phi_2$ are \textit{strongly mutually synchronous
regularity preserving} if for every synchronously regular relation 
$R \subseteq X^+ \times X^+$, the relation
$$\lbrace (u \phi_1, v \phi_2) \mid (u, v) \in R \cap (A_1 \times A_2) \
  \rbrace \subseteq Y^+ \times Y^+$$
is synchronously regular.

A \textit{sliding window inverse} for a
function $\phi : A \to Y^+$ is a
quadruple $(n, f, g, h)$ consisting of a positive integer $n$ and three   
functions
$$f : \OpPref_n(A \phi \cap Y^{> n}) \to X^*,$$
$$g : \OpInt_n(A \phi) \to X^* \text{ and}$$
$$h : \OpSuff_n(A \phi \cap Y^{> n}) \to X^*,$$ with
the property that for any $w \in A$ and $y_1 \dots y_k \in Y^+$ with
$y_1, \dots, y_k \in Y^1$ and $k > n$ such that
$w \phi = y_1 \dots y_k$, we have
\begin{align*}
w &= (y_1 \dots y_n) f \ (y_2 \dots y_{n+1}) g \ (y_3 \dots y_{n+2}) g \ 
 \dots \\
 & \hspace{18ex} \dots \ (y_{k-n} \dots y_{k-1}) g \ (y_{k-n+1} \dots y_k) h.
\end{align*}
If $f$, $g$ and $h$ are functions with domains containing those given above,
we shall say that $(n, f, g, h)$ is a sliding window inverse for $\phi$ if
the restrictions of $f$, $g$ and $h$ to the appropriate domains have the
given properties.

\begin{lemma}{\cite[Lemma~3.8]{KambitesAutoCat}} \label{lem:regpreserve}
Let $X$ and $Y$ be finite graphs, $A \subseteq X^*$ and $\phi : A \to Y^*$
be a function. If $A \phi \subseteq Y^*$ is regular, and $\phi$
has a sliding window inverse, then $\phi$ is strongly regularity preserving.
\end{lemma}

Let $\psi : \BN \to \BN$ be a function. We say that sliding window 
inverses $(n_1, f_1, g_1, h_1)$ for $\phi_1$ and $(n_2, f_2, g_2, h_2)$
for $\phi_2$ are \textit{synchronised by $\psi$} if
\begin{itemize}
\item[(i)] $n_1 = n_2$;  
\item[(ii)] for every $i = 1, 2$ and every
$w \in \OpPref_{n_1}(A_i \phi_i \cap Y^{> n_1})$, we have
$|w f_i| = 0 \psi$;
\item[(iii)] for every $i=1, 2$ and every $xyz \in A_i \phi_i$ with
             $x \in Y^m$, $m \geq 1$, $y \in Y^{n_1}$ and $z \in Y^+$
             we have $|y g_i| = m \psi$; and
\item[(iv)]  for every $i = 1, 2$ and every $xy \in A_i \phi_i$ with   
             $x \in Y^m$, $m \geq 1$ and $y \in Y^{n_1}$ we have
             $|y h_i| \leq m \psi$.
\end{itemize}
We say that two sliding window inverses are \textit{synchronised} if they
are synchronised by some function $\psi : \BN \to \BN$.  

\begin{lemma}{\cite[Lemma~3.9]{KambitesAutoCat}} \label{lem:syncratpreserve}
Let $X$ and $Y$ be finite graphs, $A_1, A_2 \subseteq X^+$ and
$\phi_1 : A_1 \to Y^+$
and $\phi_2 : A_2 \to Y^+$ be injective functions. If $A_1 \phi_1$ and
$A_2 \phi_2$ are
regular and $\phi_1$ and $\phi_2$ have synchronised sliding window inverses,
then $\phi_1$ and $\phi_2$ are strongly mutually synchronous regularity
preserving.  
\end{lemma}  

A \textit{choice of representatives} for a semigroupoid $S$ is a triple
$(X, K, \rho)$, of a graph $X$ with $X^0 = S^0$, a (surjective)
semigroupoid morphism $\rho : X^+ \to S$ of the free semigroupoid $X^+$
onto $S$, and a language $K \subseteq X^+$ such that $K \rho = S$. The
choice of representatives is called \textit{finitely generated} if $X$ has
finitely many edges. Clearly, a semigroupoid has a finitely generated
choice of representatives if and only if it is finitely generated. The
choice of representatives is called a \textit{cross-section} if the
restriction of $\rho$ to $K$ is bijective, that is, if $K$ contains a
unique representative for every arrow in $S$. The choice of
representatives is called \textit{prefix-closed} if $K$ is closed under
the taking of non-empty prefixes.

An \textit{automatic structure} for a semigroupoid $S$ is a
finitely generated choice of representatives $(X, K, \rho)$ with the
property that for every edge or empty path $a \in X^0 \cup X^1$, the
binary relation $$K_a = \lbrace (u, v) \in K \times K \mid u \omega = a
\alpha, (ua) \rho = v \rho \rbrace$$ is synchronously regular.
Equivalently \cite[Proposition~4.1]{KambitesAutoCat}, $(X, K, \rho)$ is
an automatic structure exactly if $K_a$ is synchronously regular for
every $a \in X^1$, and the union
$$K_= = \lbrace (u,v) \in K \times K \mid u \rho = v \rho  \rbrace = \bigcup_{a \in X^0} K_a$$
is synchronously regular.

A \textit{prefix-automatic structure} for a semigroupoid $S$ is an automatic
structure $(X, K, \rho)$ with the additional property that the binary relation
$$K_=' = \lbrace (u, v) \in K \times \OpPref(K) \mid u \rho = v \rho \rbrace$$
is synchronously regular. A semigroupoid admits a prefix-automatic structure
if and only if it admits an automatic structure which is prefix-closed
\cite[Corollary~4.6]{KambitesAutoCat}.

Now let $S$ be a semigroupoid and $0$ be a new symbol not in $S^1$. The
\textit{consolidation} of $S$ is the semigroup with set of elements
$T = S^1 \cup \lbrace 0 \rbrace$, and multiplication defined by
$$ s t = \begin{cases}
                     \text{the $S$-product } st & \text{ if } \
                          s, t \in S \text{ and } s \omega = t \alpha \\
                     0                            & \text{ otherwise} \\
       \end{cases}$$
for all $s, t \in S$. We shall need the following key results from \cite{KambitesAutoCat}.
\begin{theorem}{\cite[Theorem~4.3]{KambitesAutoCat}}\label{thm:cons_iff_sgpoid}
Let $T$ be the consolidation of a semigroupoid $S$. Then $S$ is automatic
[prefix-automatic] if and only if $T$ is automatic [respectively, prefix-automatic].
\end{theorem}

\begin{theorem}{\cite[Theorem~5.6]{KambitesAutoCat}}\label{thm:auto_cofinite_subsgpoid}
Let $S$ be a semigroupoid and $U$ a cofinite subsemigroupoid of $S$.
Then $S$ is automatic [prefix-automatic] if and only if $U$ is automatic
[respectively, prefix-automatic].
\end{theorem}

\section{Some Technical Results}\label{sec:technical}

We shall need the following technical results in the remaining sections.

\begin{proposition}\label{prop:syncrat_changelastletter}
Let $X$ be a finite graph and $R \subseteq X^+ \times X^+$ be synchronously
regular. Let $m \geq 0$ and let $f : X^m \to X^*$ be such that
$x f \alpha = x \alpha$ for all $x \in X^m$. Then the binary relation
$$R' = \lbrace (u, v (xf)) \mid (u, vx) \in R, u, v \in X^*, \
 x \in X^m \rbrace$$
is synchronously regular.
\end{proposition}
\begin{proof}
First, notice that $R'$ is the composition of $R$ with the relation
$$T = \lbrace \left( vx, v(xf) \right) \mid v \in X^* \rbrace.$$
By Lemma~\ref{lem:syncrat_prop}, it will suffice to show that $T$ is
synchronously rational. But clearly $T \delta_X$ is the concatenation
of $D \delta_X$ with $E \delta_X$ where
$D$ is the diagonal relation on $X$ and $E$ is the finite relation
$\lbrace (x, xf) \mid x \in X^m \rbrace$. Both of these are regular,
and regular path languages are closed under concatenation, so it follows
that $T \delta_X$ is regular, and so $T$ is synchronously regular as
required.
\end{proof}

\begin{proposition}\label{prop:no_duplicate_letters}
Let $S$ be a semigroupoid with an automatic cross-section [prefix-closed automatic
structure, prefix-automatic cross-section], and let $T$ be a finite subset
of $S$. Then $S$ has an
automatic cross-section [prefix-closed automatic structure, prefix-automatic cross-section]
$(Y, L, \sigma)$ such
that the restriction of $\sigma$ to $Y$ is injective and has image containing $T$.
\end{proposition}
\begin{proof}
Let $(X, K, \rho)$ be an automatic cross-section
 [prefix-closed automatic structure, prefix-automatic cross-section]
for $S$. Choose a subgraph $Y$ of $X$ such that $Y^1 \rho = X^1 \rho$ and
the restriction of $\rho$ to $Y^1$ is injective. For each element
$t \in T$ which does not already have a representative in $X$, adjoin
a new edge $y_t$ to $Y$ with $y_t \alpha = t \alpha$ and $y_t \omega = t \omega$,
to obtain a new finite graph $Z$.

For each $x \in X^1$, let
$x \sigma$ be the unique element $y \in Z$ such that $y \rho = x \rho$.
Extend $\sigma$ to a semigroupoid morphism $\sigma' : X^+ \to Z^+$.
Let $\rho' : Z^+ \to S$ be defined by
$$x \rho' = \begin{cases} x \rho &\text{ if } x \in X \\
                          t &\text{ if } x = y_t.
            \end{cases}
$$
It is a routine exercise to verify that $(Z, K \sigma', \rho')$ is an
automatic cross-section [prefix-closed automatic structure, prefix-automatic
cross-section] for $S$ and that the restriction of $\rho'$ to $Z$ is
injective and has image containing $T$.
\end{proof}

\begin{lemma}\label{lem:limitedwork1}
Let $(X, K, \sigma)$ be a regular choice of representatives for a
semigroupoid $S$, and let $L$ be a cofinite subset of $K$. If
$K_= \cap (L \times L)$ is synchronously regular and $K_c \cap (L \times K)$
is synchronously regular for every $c \in X^1$, then $(X, K, \sigma)$ is
an automatic structure for $S$.
\end{lemma}
\begin{proof}
First, we claim that for any $w \in K$, the language
$$w \sigma \sigma^{-1} \cap K = \lbrace x \in K \mid x \sigma = w \sigma \rbrace$$
is regular.
If $L$ contains no representatives for $w \sigma$, then
$w \sigma \sigma^{-1} \cap K \subseteq K \setminus L$ must be finite and
hence regular. Otherwise, let $u$ be a representative in
$L$ for $w \sigma$. Now we have
$$w \sigma \sigma^{-1} \cap K = (w \sigma \sigma^{-1} \cap (K \setminus L)) \cup ((K_= \cap (L \times L)) \cap (\lbrace u \rbrace \times L)) \pi_2$$
where $w \sigma \sigma^{-1} \cap K \setminus L \subseteq K \setminus L$ is finite,
$K_= \cap (L \times L)$ is synchronously regular by assumption, and $L$ is
regular. It follows by Lemma~\ref{lem:syncrat_prop}
that $w \sigma \sigma^{-1} \cap K$ is regular.
Now we see that
\begin{align*}
K_= = (K_= \cap (L \times L)) \cup \
        &\left( \bigcup_{w \in K \setminus L} \lbrace w \rbrace \times (w \sigma \sigma^{-1} \cap K) \right) \cup \\
        &\left( \bigcup_{w \in K \setminus L} (w \sigma \sigma^{-1} \cap K) \times \lbrace w \rbrace \right)
\end{align*}
is synchronously regular by Lemma~\ref{lem:syncrat_prop}(iii) and (iv).
Similarly, for any $c \in C^1$ we see that
$$K_c = (K_c \cap (L \times K)) \cup \
        \bigcup_{w \in K \setminus L} \lbrace w \rbrace \times ((wc) \sigma \sigma^{-1} \cap K)$$
is synchronously regular by Lemma~\ref{lem:syncrat_prop}(iii) and (iv).
Thus, $(X, K, \sigma)$ is an automatic structure for $S$, as required.
\end{proof}

\begin{lemma}\label{lem:limitedwork2}
Let $(X, L, \sigma)$ be an automatic structure for a semigroupoid $S$,
and suppose $L \subseteq K \subseteq X^+$ with $K \setminus L$ finite. Then
$(X, K, \sigma)$ is an automatic structure for $S$.
\end{lemma}
\begin{proof}
First, note that $K_= \cap (L \times L) = L_=$ is synchronously regular
and $K_c \cap (L \times L) = L_c$ is synchronously regular for every $c \in X^1$.
We can deduce as in the proof of Lemma~\ref{lem:limitedwork1} that
$w \sigma \sigma^{-1} \cap K$ is regular for every $w \in K$, and that
$K_=$ is synchronously regular.

Now let $c \in X^1$. For any $w \in K$, choose $u \in L$ such that
$w \sigma = u \sigma$, and consider the set
\begin{align*}
&\lbrace x \in K \mid x \omega = c \alpha, (xc) \sigma = w \sigma \rbrace \\
= &\lbrace x \in K \setminus L \mid x \omega = c \alpha, (xc) \sigma = w \sigma \rbrace
\cup \lbrace x \in L \mid x \omega = c \alpha, (x c) \sigma = w \sigma \rbrace \\
= &\lbrace x \in K \setminus L \mid x \omega = c \alpha, (x c) \sigma = w \sigma \rbrace
\cup ((K_c \cap (L \times L)) \cap (L \times \lbrace u \rbrace)) \pi_1.
\end{align*}
Certainly the set
$\lbrace x \in K \setminus L \mid x \omega = c \alpha, (x c) \sigma = w \sigma \rbrace$
is finite, and we know that $K_c \cap (L \times L)$ is synchronously
regular and $L$ is regular. It follows by Lemma~\ref{lem:syncrat_prop}(ii)
(iii) and (iv) that
$\lbrace x \in K \mid x \omega = c \alpha, (x c) \sigma = w \sigma \rbrace$ is regular.

It follows that
\begin{align*}
K_c = (K_c \cap (L \times L)) \ \cup \
        &\bigcup_{w \in K \setminus L} \lbrace w \rbrace \times ((wc) \sigma \sigma^{-1} \cap K) \ \cup \\
        &\bigcup_{w \in K \setminus L} \lbrace x \in K \mid x \omega = c \alpha, (x c) \sigma = w \sigma \rbrace \times \lbrace w \rbrace
\end{align*}
is synchronously regular by Lemma~\ref{lem:syncrat_prop}(iii) and (iv).
Thus, $(X, K, \sigma)$ is an automatic structure for $S$, as required.
\end{proof}

\section{Rees Matrices and Automaticity}\label{sec:autorees}

Rees matrix constructions over semigroupoids were introduced by
Lawson \cite{Lawson00}, who used them to construct a class of locally adequate
abundant semigroups. They represent an alternative formulation of certain special
cases of the \textit{blocked Rees matrix semigroup} constructions introduced by
Fountain \cite{Fountain82} and subsequently employed by Armstrong
\cite{Armstrong86}. In this section, we extend some results of
\cite{Descalco01} concerning Rees matrix constructions
over semigroups to cover similar cosntructions over semigroupoids and small categories. In some cases, our new results are stronger than previously known results, even when specialised to the case of semigroups.

Let $S$ be a non-empty, isolation-free semigroupoid, and $I$ and
$\Lambda$ be indexing sets. Let $F : I \to S^1 \alpha \subseteq S^0$
and
$G : \Lambda \to S^1 \omega \subseteq S^0$
be surjective functions.
Let $0$ be a new symbol not in $S$, and let
$P$ be a $\Lambda \times I$ matrix with entries drawn from $S \cup \lbrace 0 \rbrace$,
with the property that $P_{\lambda i} \omega = iF$ and $P_{\lambda i} \alpha = \lambda G$
for all $i \in I, \lambda \in \Lambda$ such that $P_{\lambda i} \neq 0$.
The \textit{Rees matrix semigroup with zero $M = M^0(S; F(I), G(\Lambda); P)$} is the semigroup with
set of elements
$$M = \lbrace 0 \rbrace \cup \lbrace (i, x, \lambda) \in I \times S \times \Lambda \mid i F = x \alpha, x \omega = \lambda G \rbrace$$
and multiplication given by
$$(i, x, \lambda)(j, y, \mu) = \begin{cases}
                                   (i, x P_{\lambda j} y, \mu) & \text{ if } P_{\lambda j} \neq 0 \\
                                   0                           & \text{ otherwise}
                               \end{cases}
$$
for all $(i, x, \lambda), (j, y, \mu) \in M \setminus \lbrace 0 \rbrace$, and
$0 m = m 0 = 0$
for all $m \in M$.

Note that in the expression ``$M^0(S; F(I), G(\Lambda); P)$'', the use of the
notation ``$F(I)$'' and ``$G(\Lambda)$'' is purely symbolic. It is intended
to remind the reader of the relationship between $F$ and $I$ and between $G$
and $\Lambda$.

We call $P$ the \textit{sandwich matrix} of
the construction. If $P$ contains no zero entries, then $M \setminus \lbrace 0 \rbrace$
is a subsemigroup of $M$, which we call the \textit{Rees matrix semigroup
(without zero) $M(S; F(I), G(\Lambda); P)$}.

We remark briefly upon the choice of codomains for the indexing functions
$F$ and $G$. This restriction does not limit the
range of semigroups which appear as Rees matrix semigroups. Indeed, if
$I$, $\Lambda$, $F$, $G$, $S$ and $P$, do not satisfy this requirement
but otherwise satisfy the requirements for the Rees matrix construction,
then one can instead perform the construction with subsets of $I$ and
$\Lambda$, a subsemigroupoid of $S$ and the corresponding restrictions
of $F$ and $G$ and submatrix of $P$, to obtain the same Rees matrix
semigroup as would be obtained by using the construction on the more
general semigroupoid. Similarly, given a semigroupoid $S$ with isolated objects, one
can remove those objects to obtain an isolation-free semigroupoid before
using a Rees matrix construction. The purpose of the restrictions is to
allow the following straightforward proposition.

\begin{proposition}\label{prop:all_sgpoid_used}
Let
$$M = M^0(S; F(I), G(\Lambda); P)$$
or
$$M = M(S; F(I), G(\Lambda); P)$$
be a Rees matrix semigroup with or without zero over a semigroupoid $S$.
Then
\begin{itemize}
\item for every $i \in I$ there exist $s \in S^1$ and $\lambda \in \Lambda$
with $(i, s, \lambda) \in M$;
\item for every $s \in S^1$ there exist $i \in I$ and $\lambda \in \Lambda$
with $(i, s, \lambda) \in M$;
\item for every $\lambda \in \Lambda$ there exist $i \in I$ and $s \in S^1$
with $(i, s, \lambda) \in M$; and
\item for every $v \in S^0$ there exists $i \in I$ with $i F = v$ or
$\lambda \in \Lambda$ with $\lambda G = v$.
\end{itemize}
\end{proposition}

The following theorem is an amalgamation of results which can be found in
\cite{Kambites03}; the same results appear, in slightly less generality, in
\cite{Kambites05}. For a precise definition of what it means for a
semigroupoid to be finitely presentable, the reader is directed to
\cite{Kambites05}.

\begin{theorem}\label{thm:fgp_rees}
Let
$$M = M(S; F(I), G(\Lambda); P)$$ or
$$M = M^0(S; F(I), G(\Lambda); P)$$ be a Rees matrix
semigroup (with or without zero) over an isolation-free semigroupoid $S$. Then $M$ is
finitely generated [finitely presentable] if and only if
\begin{itemize}
\item[(i)] the indexing sets $I$ and $\Lambda$ are finite;
\item[(ii)] the semigroupoid $S$ is finitely generated [respectively, finitely presentable]; and
\item[(iii)] $S P' S$ is a cofinite subsemigroupoid of $S$,
where $P'$ is the set of non-zero entries in the sandwich matrix $P$.
\end{itemize}
\end{theorem}

In order for a semigroup or semigroupoid to be automatic it is, of course, necessary
that it be finitely generated. Theorem~\ref{thm:fgp_rees} gives necessary and sufficient
conditions for a Rees matrix semigroup over a finitely generated semigroupoid to be finitely generated. We begin by showing that the same conditions suffice to ensure that a Rees matrix semigroup over an automatic small category is automatic. We shall then extend this result to cover Rees matrix semigroups over general automatic semigroupoids.

\begin{theorem}\label{thm:auto_cat_imp_rees}
Let $M = M^0(S; F(I), G(\Lambda); P)$ be a Rees matrix semigroup with zero over a
small category $S$. Suppose $S$ is automatic (as a semigroupoid) and $M$ is finitely
generated. Then $M$ is automatic.
\end{theorem}
\begin{proof}
Let $P'$ be the set of non-zero entries in the sandwich matrix $P$, and let
$U = S P' S$. Since $M$ is finitely generated and $S$ is a small category, we
deduce, by Theorem~\ref{thm:fgp_rees}, that
$I$ and $\Lambda$ are finite and that $U$ is a cofinite subsemigroupoid of $S$. It follows by
Theorem~\ref{thm:auto_cofinite_subsgpoid} that $U$
is automatic, by \cite[Corollary~4.6]{KambitesAutoCat} that $U$ has
an automatic cross-section, and then by Proposition~\ref{prop:no_duplicate_letters}
that $U$ has an automatic cross-section $(X, K, \sigma)$ with the property
that the restriction of $\sigma$ to $X$ is injective.

By the definition
of $U$, for each $y \in X$ we can choose $s_y, t_y \in S$, $i_y \in I$
and $\lambda_y \in \Lambda$ such that $y \sigma = s_y P_{\lambda_y i_y} t_y$.
Now let
\begin{equation*}\label{eqn:hdefine}
H = \lbrace t_y, s_y \mid y \in X \rbrace \cup \lbrace 1_s \mid s \in S^0 \rbrace \subseteq S.
\end{equation*}
We define new alphabets
$$A = \lbrace a_{i g h \lambda} \mid i \in I, g, h \in H, \lambda \in \Lambda, i F = g \alpha, g \omega = h \alpha, h \omega = \lambda G \rbrace$$
and
$$B = \lbrace b_{i s \lambda} \mid i \in I, s \in S \setminus U, \lambda \in \Lambda, i F = s \alpha, s \omega = \lambda G \rbrace.$$
Clearly, $A$ and $B$ are finite. Let $z$ be a new symbol not in $A$ or $B$, which will represent the zero element
$0 \in M$. We define a morphism of semigroups
$$\rho : (A \cup B \cup \lbrace z \rbrace)^+ \to M$$
by $(a_{i g h \lambda}) \rho = (i, gh, \lambda)$, $(b_{i s \lambda}) \rho = (i, s, \lambda)$
and $z \rho = 0$.

For every $i \in I$ and $\lambda \in \Lambda$, let $X_{i \lambda}$ be
the set
$$X_{i \lambda} = \lbrace x \in X^+ \mid i F = x \alpha, x \omega = \lambda G \rbrace$$
of words in $X^+$ starting at $i F$ and finishing at $\lambda G$,
and define a function
$\phi_{i \lambda} : X_{i \lambda} \to A^+$
by
$$(w_1 w_2 \dots w_n) \phi_{i \lambda} = a_{i,1_{i F},s_{w_1},\lambda_{w_1}} a_{i_{w_1},t_{w_1},s_{w_2},\lambda_{w_2}} \dots a_{i_{w_n}, t_{w_n}, 1_{\lambda G}, \lambda}$$
for all $w_1, w_2, \dots, w_n \in X$ with $w_1 \dots w_n \in X_{i \lambda}$.
It follows from the fact that the restriction of $\sigma$ to $X$ is injective,
that each function $\phi_{i \lambda}$ is injective. Indeed, if
$(w_1 \dots w_n) \phi_{i \lambda} = (w'_1 \dots w'_{n'}) \phi_{i \lambda}$
then from the definition of $\phi_{i \lambda}$ we have
\begin{align*}
&a_{i,1_{i F},s_{w_1},\lambda_{w_1}} a_{i_{w_1},t_{w_1},s_{w_2},\lambda_{w_2}} \dots a_{i_{w_n}, t_{w_n}, 1_{\lambda G}, \lambda}
= (w_1 \dots w_n) \phi_{i \lambda} \\
= &(w'_1 \dots w'_{n'}) \phi_{i \lambda}
= a_{i,1_{i F},s_{w'_1},\lambda_{w'_1}} a_{i_{w'_1},t_{w'_1},s_{w'_2},\lambda_{w'_2}} \dots a_{i_{w'_{n'}}, t_{w'_{n'}}, 1_{\lambda G}, \lambda}
\end{align*}
so clearly $n = n'$ and for $1 \leq k \leq n$ we have
$s_{w_k} = s_{w'_k}$, $i_{w_k} = i_{w'_k}$, $\lambda_{w_k} = \lambda_{w'_k}$
and $t_{w_k} = t_{w'_k}$. Now each
$$w_k \sigma = s_{w_k} P_{\lambda_{w_k} i_{w_k}} t_{w_k} = s_{w'_k} P_{\lambda_{w'_k} i_{w'_k}} t_{w'_k} = w'_k \sigma.$$
But $\sigma$ is by assumption injective on $X$, so it follows that
each $w_k = w'_k$.

Note also that
\begin{equation}\label{eqn:autoreesa}
(i, w \sigma, \lambda) = w \phi_{i \lambda} \rho
\end{equation}
for all paths $w \in X_{i \lambda}$.

Our aim is to show that each function of the form $\phi_{i \lambda}$ is
strongly regularity preserving, and that each pair of such functions is
strongly mutually synchronous regularity preserving. We shall do so by
showing that each such function has a regular image and that the functions
have sliding window inverses which are pairwise synchronised, and then
invoking Lemmas~\ref{lem:regpreserve} and \ref{lem:syncratpreserve}.

We claim first that the image $X_{i \lambda} \phi_{i \lambda}$ of each
$\phi_{i \lambda}$ is regular. We say that a two-letter word
$a_{i_1, g_1, h_1, \lambda_1} a_{i_2, g_2, h_2, \lambda_2} \in A^2$
is \textit{compatible} if there exists $y \in X$ such that
$s_y = h_1$, $\lambda_y = \lambda_1$, $i_y = i_2$ and $t_y = g_2$, and
\textit{incompatible} otherwise.
Then clearly, $X_{i \lambda} \phi_{i \lambda}$ is the set of words in $A^{\geq 2}$
which begin with a letter of the form $a_{i,1_{i F}, h, \mu}$ for some
$h \in H$ and $\mu \in \lambda$, end with a letter of the form
$a_{j, g, 1_{\lambda G}, \lambda}$ for some $j \in I$ and $g \in H$, and
contain no incompatible factors of length $2$. Thus,
$X_{i \lambda} \phi_{i \lambda} = (A_1 A^* A_2) \setminus (A^* A_3 A^*)$
where
$$A_1 = \lbrace a_{i, 1_{i F}, h, \mu} \in A \mid h \in H, \mu \in \Lambda, i F = h \alpha, h \omega = \mu G \rbrace,$$
$$A_2 = \lbrace a_{j, g, 1_{\lambda G}, \lambda} \in A \mid j \in I, g \in H, j F = g \alpha, g \omega = \lambda G \rbrace, \text{ and}$$
$$A_3 = \lbrace ab \in A^2 \mid ab \text{ is not compatible} \rbrace$$
are all finite.
It follows that $X_{i \lambda} \phi_{i \lambda}$ is regular, as required.

Next, we claim that the $\phi_{i \lambda}$ have
pairwise synchronised sliding window inverses.
Fix $i \in I$ and $\lambda \in \Lambda$. We shall define a function
$$f : \lbrace ab \in A^2 \mid ab \text{ is compatible} \rbrace \to X.$$
By our observations above, the domain of $f$ will then include all
two-letter factors of words in $X_{i \lambda} \phi_{i \lambda}$.

Suppose $a_{i_1 g_1 h_1 \lambda_1} a_{i_2 g_2 h_2 \lambda_2}$ is compatible.
Then by definition, there is some $y \in X$ with $s_y = h_1$, $\lambda_y = \lambda_1$, $i_y = i_2$
and $t_y = g_2$. Furthermore, we have $y \sigma = h_1 P_{\lambda_1 i_2} g_2$, and
$\sigma$ is injective on $X$, so $y$ is the unique letter in $X$ with this property.
Thus, we can well-define $f$ by the rule
$$(a_{i_1 g_1 h_1 \lambda_1} a_{i_2 g_2 h_2 \lambda}) f = y.$$

We claim that $(2, f, f, f)$ is a sliding window inverse
for $\phi_{i \lambda}$.
To show this, suppose $w = w_1 w_2 \dots w_n \in X_{i \lambda}$ with
$w_1, \dots, w_n \in X$. Then by the definition of $\phi_{i \lambda}$ we have
$$(w_1 w_2 \dots w_n) \phi_{i \lambda} = a_{i,1_{i F},s_{w_1},\lambda_{w_1}} a_{i_{w_1},t_{w_1},s_{w_2},\lambda_{w_2}} \dots a_{i_{w_n}, t_{w_n}, 1_{\lambda G}, \lambda}.$$
Let
$$y_1 = a_{i,1_{i F},s_{w_1},\lambda_{w_1}},
\hspace{1em} y_2 = a_{i_{w_1},t_{w_1},s_{w_2},\lambda_{w_2}},
\hspace{1em} \dots,
\hspace{1em} y_{n+1} = a_{i_{w_n}, t_{w_n}, 1_{\lambda G}, \lambda}$$
so that
$$(w_1 w_2 \dots w_n) \phi_{i \lambda} = y_1 y_2 \dots y_{n+1}.$$
Now from the definition of $f$, we see that $(y_i y_{i+1}) f = w_i$ for
each $1 \leq i \leq n$, so we have
$$(y_1 y_2) f (y_2 y_3) f (y_3 y_4) f \dots (y_{n-1} y_{n}) f (y_{n} y_{n+1}) f = w_1 \dots w_n$$
as required. It follows by Lemma~\ref{lem:regpreserve} that each $\phi_{i \lambda}$ is
strongly regularity preserving.

Furthermore, if we define $\psi : \BN \to \BN$ to be the constant function
given by $n \psi = 1$, then $m \psi = 1 = |w f|$
for any $w \in \OpFact_2(X_{i \lambda} \phi_{i \lambda})$
which occurs starting in position $m$ (numbered from zero).
Since the definition of $\psi$ is independent of the choice of $i$ and $\lambda$, it follows
that $\psi$ synchronises the sliding window inverses constructed above for any pair of functions
of the form $\phi_{i \lambda}$. We have already observed that the
$\phi_{i \lambda}$ are injective and have regular images, so it follows by
Lemma~\ref{lem:syncratpreserve} that
any pair of such functions are strongly mutually synchronous regularity preserving.

We now define a language $L \subseteq (A \cup B \cup \lbrace z \rbrace)^+$ by
$$L = B \cup \lbrace z \rbrace \cup \bigcup_{i \in I, \lambda \in \Lambda} (K \cap X_{i \lambda}) \phi_{i \lambda}$$
We claim that $(A \cup B \cup \lbrace z \rbrace, L, \rho)$ is an automatic
structure (indeed, an automatic cross-section) for $M$.

We show first that $L$ contains a unique representative for every element of $M$.
Certainly $z \in L$ is the unique word representing the zero element $0 \in M$.
Now consider a non-zero element $(i, s, \lambda) \in M$. If $s \in S \setminus U$ then
$b_{i s \lambda} \in B$ is clearly the
unique word in $L$ representing $(i, s, \lambda)$, as required.
Otherwise, we have $s \in U$. Now $(X, K, \sigma)$ is a choice of representatives
for $U$, so
we must have
$s = w \sigma$
for some path $w \in K$. Furthermore, $w \alpha = s \alpha = i F$ and
$w \omega = s \omega = \lambda G$, so that $w \in X_{i \lambda}$.
Now by \eqref{eqn:autoreesa} we have
$$(i, s, \lambda) = (i, w \sigma, \lambda) = w \phi_{i \lambda} \rho$$
where $w \phi_{i \lambda} \in (K \cap X_{i \lambda}) \phi_{i \lambda} \subseteq L$.
Furthermore, if $i' \in I$ and $\lambda' \in \Lambda$ and $v \in K \cap X_{i' \lambda'}$ are such that
$v \phi_{i' \lambda'} = (i, s, \lambda)$, then by \eqref{eqn:autoreesa} we must have
$(i', v \sigma, j') = (i, w \sigma, j)$
from which it follows that $i = i'$, $j = j'$ and, since $K$ contains a unique
representative for every element in $U$, $v = w$. Hence $w \phi_{i \lambda}$ is the
unique representative for $(i, s, \lambda)$ in $L$. Thus,
$(A \cup B \cup \lbrace z \rbrace, L, \rho)$ is a cross-section for $M$.

Next, we show that $L$ is a regular language.
Clearly, $B$ and $\lbrace z \rbrace$ are finite and hence regular. Furthermore, $K$ is
regular and each $\phi_{i \lambda}$ is strongly regularity preserving, so each
$(K \cap X_{i \lambda}) \phi_{i \lambda}$ is regular. It follows that
$L$ is a finite union of regular languages, and hence is regular.

Since $(A \cup B \cup \lbrace z \rbrace, L, \rho)$ is a cross-section for
$M$, it follows by Lemma~\ref{lem:syncrat_prop}(v) that the binary relation
$$L_= = \lbrace (u, v) \in L \times L \mid u \rho = v \rho \rbrace = \lbrace (u, u) \mid u \in L \rbrace$$
is synchronously regular.

It remains to show that $L_a$ is synchronously regular for all
$a \in A \cup B \cup \lbrace z \rbrace$.
First, since $z$ is the unique representative in $L$ for the zero element
$0 \in M$, we have
$L_z = \lbrace (w, z) \mid w \in L \rbrace = L \times \lbrace z \rbrace$
which is synchronously regular by Lemma~\ref{lem:syncrat_prop}(iii).

Now suppose $a \in A \cup B$,
with $a \rho = (i_a, s_a, \lambda_a)$. By Lemma~\ref{lem:limitedwork1},
it will suffice to show that
$$L_a \cap ((A^+ \cap L) \times L) = L_a \cap (A^+ \times L)$$ is
synchronously regular. We write $L_a \cap (A^+ \times L)$ as a union
\begin{align*}
L_a \cap (A^+ \times L) = &(L_a \cap (A^+ \times B)) \\
       &\cup (L_a \cap (A^+ \times \lbrace z \rbrace)) \\
       &\cup (L_a \cap (A^+ \times A^+)),
\end{align*}
and show that each component in this union is synchronously regular.

Suppose $(u, v) \in L_a \cap (A^+ \times L)$, so that
$(u \rho) (a \rho) = v \rho$.
Clearly $v \notin B$, since letters in $B$ represent indecomposable elements
of $M$, so $L_a \cap (A^+ \times B)$ is empty and hence synchronously regular.

If $v = z$, then we have $(u \rho) (a \rho) = v \rho = z \rho = 0$. Since
$u \in A^+ \cap L$ cannot represent the zero element, this can happen only if $P_{u \rho \pi_3, a \rho \pi_1} = 0$. Thus,
\begin{align*}
L_a \cap (A^+ \times \lbrace z \rbrace) &= \lbrace (u, z) \mid u \in A^+, P_{u \rho \pi_3, a \rho \pi_1} = 0 \rbrace \\
                                       &= \left(\bigcup_{b \in A, P_{(b \rho \pi_3), a \rho \pi_1} = 0} (A^* b \cap L) \right) \times \lbrace z \rbrace
\end{align*}
is a product of finite unions of regular languages, and by
Lemma~\ref{lem:syncrat_prop}(iii) is synchronously regular.

It remains only to show that $L_a \cap (A^+ \times A^+)$ is synchronously
regular. Certainly we have
$$L_a \cap (A^+ \times A^+) = \bigcup_{i \in I, \lambda \in \Lambda} L_{a, i, \lambda}$$
where each
$$L_{a,i,\lambda} = \lbrace (u, v) \in L_a \mid u, v \in A^+, u \rho \pi_1 = v \rho \pi_1 = i, u \rho \pi_3 = \lambda \rbrace.$$
If $P_{\lambda, (a \rho \pi_1)} = 0$ then $L_{a, i, \lambda}$ is empty and hence
synchronously regular. Otherwise, let $w$ be a path in $K$ such that
$w \sigma = P_{\lambda, (a \rho \pi_1)} (a \rho \pi_2)$. It follows from
\cite[Proposition~4.2]{KambitesAutoCat} the language $K_w$
is regular. Now using \eqref{eqn:autoreesa} we have
\begin{align*}
L_{a, i, \lambda} &= \lbrace (u, v) \in (K \cap X_{i \lambda}) \phi_{i \lambda} \times (K \cap X_{i,a \rho \pi_3}) \phi_{i, a \rho \pi_3} \mid (u \rho) (a \rho) = v \rho \rbrace \\
                  &= \lbrace (u, v) \in (K \cap X_{i \lambda}) \phi_{i \lambda} \times (K \cap X_{i,a \rho \pi_3}) \phi_{i, a \rho \pi_3} \\
&\hskip 17em \mid (u \rho \pi_2) P_{\lambda, a \rho \pi_1} (a \rho \pi_2) = v \rho \pi_2 \rbrace \\
                  &= \lbrace (x \phi_{i \lambda}, y \phi_{i, a \rho \pi_3}) \mid (x, y) \in K_w \cap (X_{i \lambda} \times X_{i, a \rho \pi_3}) \rbrace.
\end{align*}

But $\phi_{i \lambda}$ and $\phi_{i, a \rho \pi_3}$ are strongly mutually
synchronous regularity preserving, so it follows that $L_{a, i, \lambda}$
is synchronously regular. Hence, $L_a$ is a finite union of synchronously regular
binary relations, and by Lemma~\ref{lem:syncrat_prop}(iv) is synchronously
regular as required.
\end{proof}

This result extends easily from small categories to semigroupoids.

\begin{theorem}\label{thm:auto_sgpoid_imp_rees}
Let $M = M^0(S; F(I), G(\Lambda); P)$ be a Rees matrix semigroup with zero over a
semigroupoid $S$. Suppose $S$ is automatic and $M$ is finitely
generated. Then $M$ is automatic.
\end{theorem}
\begin{proof}
Let $P'$ be the set of non-zero entries in the sandwich matrix $P$ and
let $U = S P' S$.
Since $M$ is finitely generated, Theorem~\ref{thm:fgp_rees} tells
us that $I$, $\Lambda$ and $S \setminus U$ are all finite.
Let
$\overline{S}$ be the category obtained from $S$ by adjoining a new
identity arrow $1_e$ at every object $e \in S^0$ (even if there is
already an identity arrow at $e$). Since $I$ and $\Lambda$ are finite,
the images $I F$ and $\Lambda G$ of $F$ and $G$ are finite. But since
$S$ is isolation-free, it follows from
Proposition~\ref{prop:all_sgpoid_used} that $S^0$ is the union of
$I F$ and $G \Lambda$, so $S$ has only finitely many objects.
It follows that $S$ is a cofinite subsemigroupoid of
$\overline{S}$ so, by Theorem~\ref{thm:auto_cofinite_subsgpoid},
$\overline{S}$ is automatic.

Define
$\overline{I} = I \cup \lbrace i_s \mid s \in S^0 \setminus I F \rbrace$
and extend $F$ to a function $\overline{F} : \overline{I} \to S^0$ by defining
$i_s \overline{F} = s$ for all $s \in S^0 \setminus I F$.
Similarly, define
$\overline{\Lambda} = \Lambda \cup \lbrace \lambda_s \mid s \in S^0 \setminus
\Lambda G \rbrace$
and extend $G$ to a function $\overline{G} : \overline{\Lambda} \to S^0$ by
defining
$\lambda_s \overline{G} = s$ for all $s \in S^0 \setminus \Lambda G$.
Extend $P$ to a $\overline{\Lambda} \times \overline{I}$ matrix
$\overline{P}$, by defining all
new entries to be zero.

Now let $\overline{M} = M^0(\overline{S}; \overline{F}(\overline{I}),
\overline{G}(\overline{\Lambda}); \overline{P})$.
We know that $\overline{I}$ and $\overline{\Lambda}$ are finite and that
$\overline{S}$ is automatic. Furthermore, the set of
non-zero entries in $\overline{P}$ is exactly the set $P'$ of non-zero
entries in $P$. Now if we let $\overline{U} = \overline{S} \ P' \ \overline{S}$, then we have
$U \subseteq \overline{U}$, from which it follows that
$\overline{S} \setminus \overline{U}$ is finite.
It follows by Theorem~\ref{thm:auto_cat_imp_rees} that $\overline{M}$
is automatic.

Now clearly every non-zero element of $M$ is also contained
in $\overline{M}$. Moreover, the multiplication in $M$ is clearly the same
as that in $\overline{M}$, so we conclude that $M$ is a subsemigroup
of
$\overline{M}$. Furthermore, the only elements of $\overline{M}$ not in $M$
are those of the form $(i, 1_e, \lambda)$ for
$i \in \overline{I}$, $\lambda \in \overline{\Lambda}$ and $e \in S^0$ with
$$i \overline{F} = 1_e \alpha = e = 1_e \omega = \lambda \overline{G}.$$
Since $S$ has finitely many objects and $\overline{I}$ and
$\overline{\Lambda}$ are
finite it follows that $M$ is a cofinite subsemigroup of $\overline{M}$.
By Theorem~\ref{thm:auto_cofinite_subsgpoid} (or \cite[Theorem~1.1]{Hoffmann02}),
it follows that $M$ is automatic.
\end{proof}

Combining Theorem~\ref{thm:auto_sgpoid_imp_rees} with
Theorem~\ref{thm:fgp_rees}
we obtain a more explicit sufficient condition, without reference to $M$ being
finitely generated.

\begin{corollary}\label{thm:auto_rees_condition}
Let $M = M^0(S; F(I), G(\Lambda); P)$ be a Rees matrix semigroup with zero over a
semigroupoid $S$. If
\begin{itemize}
\item[(i)] the indexing sets $I$ and $\Lambda$ are finite;
\item[(ii)] the semigroupoid $S$ is automatic; and
\item[(iii)] $S \setminus SP'S$ is finite, where $P'$ is the set of
non-zero entries in the sandwich matrix $P$
\end{itemize}
then $M$ is automatic.
\end{corollary}

Next, we give a sufficient condition for the underlying semigroupoid of
an automatic Rees matrix semigroup to be automatic.

Let $M = M^0(S; F(I), G(\Lambda); P)$ be a Rees matrix semigroup with zero
over a semigroupoid $S$. Let $T$ be a subset of $S$. We say that $T$ is
\textit{strongly right-ideal-generated by a row cross-section of $P$} if there exists a subset
$\Lambda' \subseteq \Lambda$ such that the restriction of $G$ to $\Lambda'$ is
bijective, and every arrow in $T$ can be written in the form $P_{\lambda i} s$
for some $\lambda \in \Lambda'$, $i \in I$ and $s \in S^1$. We say that
$T$ is \textit{(weakly) right-ideal-generated by a row cross-section of $P$} if there
exists $\Lambda'$ as above, such that every arrow in $T$ can be written
either in the form $P_{\lambda i} s$ or in the form $P_{\lambda i}$ (or both)
for some $\lambda \in \Lambda'$, $i \in I$ and (where appropriate) $s \in S^1$.

The next result says that the set of non-identity elements of $S$ being strongly
right-ideal-generated by a row cross-section of $P$ is a sufficient condition for
automaticity in $M$ to imply automaticity in $S$. We shall subsequently
strengthen the result by weakening the hypothesis, showing that it suffices
for the set of non-identity elements of $S$ to be \textit{weakly}
right-ideal-generated by a row cross-section of $P$.

\begin{theorem}\label{thm:auto_rees_mayimp_sgpoid1}
Let $M = M^0(S; F(I), G(\Lambda); P)$ be a Rees matrix semigroup with zero over
a semigroupoid $S$. Suppose $M$ is automatic, and the set of non-identity elements
of $S$ is strongly right-ideal-generated by a row cross-section of $P$. Then $S$ is automatic.
\end{theorem}
\begin{proof}
It follows from \cite[Corollary~4.6]{KambitesAutoCat} (or \cite[Corollary~5.6]{Campbell01})
that $M$ has an automatic cross-section and then from Proposition~\ref{prop:no_duplicate_letters}
that $M$ has an automatic cross-section $(A, L, \rho)$ with the property
that the restriction of $\rho$ to $A$ is injective.
Let $\Lambda' \subseteq \Lambda$ be such that the restriction of $G$
to $\Lambda'$ is bijective, and every non-identity arrow in $S$ can be written
in the
form $P_{\lambda i} s$ for some $\lambda \in \Lambda'$, $i \in I$ and $s \in S^1$.

Choose some subset $I' \subseteq I$
such that $F$ restricts to a bijection on $I'$. Let $F' : S^1 \alpha \to I'$ and $G' : S^1 \omega \to \Lambda'$
be the inverses of the restrictions of $F$ and $G$ to $I'$ and $\Lambda'$ respectively.
We define sets
$$C = \lbrace c_s \mid s \in S, (i, s, \lambda) \in A \rho \text{ for some } i \in I, \lambda \in \Lambda \rbrace$$
and
$$D = \lbrace d_{\lambda i} \mid \lambda \in \Lambda, i \in I, P_{\lambda i} \neq 0 \rbrace.$$
Clearly, $C$ and $D$ are finite.
Let $X$ be a graph with vertex set $S^0$, and edge set $C \cup D$
where $c_s \alpha = s \alpha$, $c_s \omega = s \omega$,
$d_{\lambda i} \alpha = P_{\lambda i} \alpha$ and
$d_{\lambda i} \omega = P_{\lambda i} \omega$.
Then there is a natural morphism $\sigma : X^+ \to S$, given by
$c_s \sigma = s$ and $d_{\lambda i} \sigma = P_{\lambda i}$. Note
that $w \sigma \alpha = w \alpha$ and $w \sigma \omega = w \omega$
for every $w \in X^1$ and hence for every $w \in X^+$.

Let $V$ be the language
$$V = (I' \times S \times \Lambda') \rho^{-1}$$
of words in $A^+$ which represent non-zero elements of the form
$(i, s, \lambda) \in M$ with $i \in I'$ and $\lambda \in \Lambda'$.
Define a function $\phi : V \to X^+$ by
\begin{equation*}\label{eqn:phidefine}
(a_1 a_2 \dots a_n) \phi = c_{s_1} d_{\lambda_1 i_2} c_{s_2} d_{\lambda_2 i_3} \dots d_{\lambda_{n-1} i_n} c_{s_n}
\end{equation*}
where each $a_k \rho = (i_k, s_k, \lambda_k)$. That $\phi$ is a well-defined
function into $X^+$ follows from the fact that words in $V$ represent non-zero elements
of $M$. Furthermore, if
$$(a_1 a_2 \dots a_n) \phi = c_{s_1} d_{\lambda_1 i_2} c_{s_2} d_{\lambda_2 i_3} \dots d_{\lambda_{m-1} i_m} c_{s_m} = (b_1 b_2 \dots b_{n'}) \phi$$
then clearly $n = m = n'$ and for $1 \leq k \leq n$ we have
$$a_k \rho = (i_k, s_k, \lambda_k) = b_k \rho$$
where $i_1 = s_1 \alpha F'$ and $\lambda_m = s_m \omega G'$.
But $\rho$ is by assumption injective when restricted to the alphabet
$A$, so it follows that $\phi$ is injective.

Notice also that for any $w \in V$ we have
\begin{equation}\label{eqn:autoreesb}
w \rho = (w \phi \sigma \alpha F', w \phi \sigma, w \phi \sigma \omega G').
\end{equation}

Our aim is to show that the function $\phi$ has a regular
image and a self-synchronised sliding window inverse.
We claim first that the image $V \phi$ of the function $\phi$ is regular.
We say that an ordered triple $(d_{\lambda i}, c_s, d_{\mu j}) \in D \times C \times D$
is a \textit{valid internal triple} if there exists a letter $a \in A$ with
$a \rho = (i, s, \mu)$.
We say that $(c_s, d_{\lambda i}) \in C \times D$ is a \textit{valid start pair} if
there exists $a \in A$ with $a \rho = (i', s, \lambda)$ for some
$i' \in I'$.
We say that $(d_{\lambda i}, c_s) \in D \times C$ is a \textit{valid end pair} if
there exists $a \in A$ with $a \rho = (i, s, \lambda')$ for some
$\lambda' \in \Lambda'$.
Now it is easily verified that
$$V \phi \cap X^{\geq 5} = Q_1 (CD)^* C Q_2 \setminus (C \cup D)^*  Q_3' (C \cup D)^*$$
where
$$Q_1 = \lbrace cd \in CD \mid (c, d) \text{ is a valid start pair} \rbrace,$$
$$Q_2 = \lbrace dc \in DC \mid (d, c) \text{ is a valid end pair} \rbrace, \text{ and }$$
$$Q_3' = \lbrace cde \in DCD \mid (c, d, e) \text{ is \textit{not} a valid internal triple} \rbrace$$
are all finite. It follows that $V \phi \cap X^{\geq 5}$ is regular,
and hence that $V \phi$ is regular, as required.

Next, we claim that $\phi$ has a self-synchronised sliding window inverse.
Consider a three-letter prefix $v_1 v_2 v_3$ of a word
$$(a_1 a_2 \dots a_n) \phi \in V \phi.$$ We shall define
$$(v_1 v_2 v_3) f = a_1.$$
To show that $f$ is well-defined, suppose $(b_1 b_2 \dots b_m) \phi$ also has
the prefix $v_1 v_2 v_3$. Then by the definition of $\phi$ we have
$$(a_1 a_2 \dots a_n) \phi = c_{s_1} d_{\lambda_1 i_2} c_{s_2} d_{\lambda_2 i_3} \dots d_{\lambda_{n-1} i_n} c_{s_n}$$
where each $a_k \rho = (i_k, s_k, \lambda_k)$, and similarly
$$(b_1 b_2 \dots b_m) \phi = c_{t_1} d_{\mu_1 j_2} c_{t_2} d_{\mu_2 j_3} \dots  \dots \dots d_{\mu_{m-1} j_m} c_{t_m}$$
where each $b_k \rho = (j_k, t_k, \mu_k)$.
But now we have $c_{s_1} = v_1 = c_{t_1}$ and $d_{\lambda_1 i_2} = v_2 = d_{\mu_1 j_2}$.
But by the definitions of $C$ and $D$, it follows that $s_1 = t_1$ and
$\lambda_1 = \mu_1$. Furthermore,
we must have $i_1 F = s_1 \alpha = t_1 \alpha = j_1 F$, but $i_1, j_1 \in I'$
and $F$ is injective on $I'$, so we must have $i_1 = j_1$. Thus, we obtain
$$a_1 \rho = (i_1, s_1, \lambda_1) = (j_1, t_1, \mu_1) = b_1 \rho.$$
But $\rho$ is, by assumption,
injective on the alphabet $A$, so we must have $a_1 = b_1$, as required to show
that $f$ is well-defined.

Similarly, given a word
$$(a_1 a_2 \dots a_n) \phi = c_{s_1} d_{\lambda_1 i_2} c_{s_2} d_{\lambda_2 i_3} \dots d_{\lambda_{n-1} i_n} c_{s_n}$$
we define
$$(d_{\lambda_{p-1} i_p} c_{s_p} d_{\lambda_p i_{p+1}}) g = a_p \text{ and } \
(c_{s_p} d_{\lambda_p i_{p+1}} c_{s_{p+1}}) g = \epsilon$$
where $\epsilon$ denotes the empty word in $A^*$.
A similar argument to that for $f$ shows that $g$ is well-defined.

Finally, given a three-letter suffix $v_1 v_2 v_3$ of a word
$(a_1 a_2 \dots a_n) \phi \in V \phi,$ we shall define
$(v_1 v_2 v_3) h = a_n.$
Once again, a similar argument to that for $f$
shows that $h$ is well-defined.

Now suppose $w = a_1 \dots a_n \in V$ and $w \phi = y_1 \dots y_m$.
Then by the definition of $\phi$, we have
$$y_1 \dots y_m = (a_1 a_2 \dots a_n) \phi = c_{s_1} d_{\lambda_1 i_2} c_{s_2} d_{\lambda_2 i_3} \dots d_{\lambda_{n-1} i_n} c_{s_n}$$
where each $a_k \rho = (i_k, s_k, \lambda_k)$.
But now by the definitions of $f$, $g$ and $h$, it follows that
\begin{align*}
(y_1 y_2 y_3) f (y_2 y_3 y_4) g (y_3 y_4 y_5) g \dots (y_{m-3} &y_{m-2} y_{m-1}) g (y_{m-2} y_{m-1} y_m) h \\
&= a_1 a_2 \epsilon a_3 \epsilon \dots \epsilon a_{n-1} a_n  = w
\end{align*}
We have shown that $(3, f, g, h)$ is a sliding window inverse for $\phi$,
and by Lemma~\ref{lem:regpreserve}, it follows that $\phi$ is strongly
regularity preserving.

Furthermore, if we define a function $\psi : \BN \to \BN$ by
$$n \psi = \begin{cases}
      1 &\text{ if } n = 0 \text{ or } n \text{ is odd} \\
      0 &\text{otherwise}
   \end{cases}$$
then since factors of paths in $V \phi$ begin with a letter from $C$ exactly
if they begin in an odd position, $\psi$ synchronises $(3, f, g, h)$ with
itself. We have already observed that $\phi$ is injective and has a regular
image. It follows by Lemma~\ref{lem:syncratpreserve} that $\phi$ is
strongly synchronous regularity preserving.

We now define $K = (L \cap V) \phi \subseteq X^+$, and claim that $(A, K, \sigma)$ is an
automatic structure (indeed, an automatic cross-section) for $S$.

First, we show that $\sigma$ maps $K$ bijectively onto $S$.
To this end, let $s \in S$. Then there is an element $m = (s \alpha F', s, s \omega G') \in M$, so there exists
a word $w \in L$ representing $m$. Indeed, since $s \alpha F' \in I'$ and
$s \omega G' \in \Lambda'$, we have $w \in V$, and so $w \phi \in K$. But
by \eqref{eqn:autoreesb}, $w \phi \sigma = w \rho \pi_2 = s$, so $w \phi$ is a representative in $K$ for $s$.

Furthermore, if $v \phi \in (L \cap V) \phi = K$ also represents $s$, then we must have
$v \rho = (i', s, \lambda')$ where $i' \in I'$ and $\lambda' \in \Lambda'$. But for
$(i', s, \lambda') \in M$, we must have $i' F = s \alpha = i F$ and $\lambda' G = s \omega = \lambda G$.
Since $F$ and $G$ are bijective when restricted to $I'$ and $\Lambda'$, we must
have $i = i'$ and $\lambda = \lambda'$, and so $v \rho = w \rho$. But $v, w \in L$
and $(A, L, \rho)$ is a cross-section for $M$, so we deduce that $v = w$, and hence
$v \phi = w \phi$. Thus, $(X, K, \sigma)$ is a cross-section for $S$.

We know that $L$ is a regular language, and that $\phi$
is strongly regularity preserving, so it is immediate that $K = (L \cap V) \phi$ is a regular
language. Since $(X, K, \sigma)$ is a cross-section
for $S$, it follows that the binary relation
$$K_= = \lbrace (u, v) \in K \times K \mid u \sigma = v \sigma \rbrace = \lbrace (u, u) \mid u \in K \rbrace$$
is synchronously regular by Lemma~\ref{lem:syncrat_prop}(v).

Now let $b$ be an edge in $X$. We must show that
$$K_b = \lbrace (u, v) \in K \times K \mid u \omega = b \alpha, (u b) \sigma = v \sigma \rbrace$$
is synchronously regular.

If $b \sigma$ is an identity arrow, then it follows easily from the fact
that $K_=$ is synchronously regular that $K_b$ is synchronously regular,
and we are done.

Otherwise, $b \sigma$ is not an identity arrow in $S$. Now
by assumption, we can write $b \sigma = P_{\lambda i} c$ for some
$\lambda \in \Lambda'$, $i \in I$ and $c \in S^1$. Furthermore, we must have
$b \alpha = b \sigma \alpha = (P_{\lambda i} c) \alpha = \lambda G$, so that
$\lambda = b \alpha G'$. Let $\mu = b \omega G'$. Now $M$ has an element
$(i, c, \mu)$. It follows from \cite[Proposition~4.2]{KambitesAutoCat}
that the language
$$L_{(i,c,\mu)} = \lbrace (u, v) \in L \times L \mid (u \rho) (i, c, \mu) = v \rho \rbrace$$
is synchronously regular. We claim that
$$K_b = \lbrace (u \phi, v \phi) \mid (u, v) \in L_{(i, c, \mu)} \cap (V \times V), u \rho \pi_3 = \lambda \rbrace.$$
It will then follow that
$$K_b = \lbrace (u \phi, v \phi) \mid (u, v) \in L_{(i, c, \mu)} \cap (A^* D \times A^+) \cap (V \times V) \rbrace$$
where
$D = \lbrace a \in A \mid a \rho \pi_3 = \lambda \rbrace$
is finite, so that $A^* D$ is regular by Lemma~\ref{lem:syncrat_prop}.
By Lemma~\ref{lem:syncrat_prop}(iii) and (iv) it will follow that
$L_{(i, c, \mu)} \cap (A^* D \times A^+)$ is synchronously regular.
Finally, since $\phi$ is strongly synchronous regularity preserving, we shall deduce that
$K_b$ is synchronously regular, as required.

To prove the claim, first suppose that $(a_1, a_2) \in K_b$. Then certainly
$a_1, a_2 \in K = (L \cap V) \phi$, so $a_1 = u_1 \phi$ and $a_2 = u_2 \phi$
for some $u_1, u_2 \in L \cap V$.
Now by \eqref{eqn:autoreesb} we have
$$u_k \rho = (u_k \phi \sigma \alpha F', u_k \phi \sigma, u_k \phi \sigma \omega G')$$
for $k = 1$ and $k = 2$.
Now we have
$$u_1 \rho \pi_3 = u_1 \phi \sigma \omega G' = a_1 \sigma \omega G' = a_1 \omega G' = b \alpha G' = \lambda.$$
Similarly, we obtain
$$u_1 \rho \pi_1 = u_1 \phi \sigma \alpha F' = a_1 \sigma \alpha F' = a_2 \sigma \alpha F' = u_2 \phi \sigma \alpha F' = u_2 \rho \pi_1$$
and
$$u_2 \rho \pi_3 = u_2 \phi \sigma \omega G' = a_2 \sigma \omega G' = b \sigma \omega G' = \mu = (i, c, \mu) \pi_3.$$
Finally, we have
$$u_2 \rho \pi_2 = u_2 \phi \sigma = a_2 \sigma = (a_1 \sigma) (b \sigma) = (a \sigma) P_{\lambda i} c = ((u_1 \rho) (i, c, \mu)) \pi_2.$$
We have shown that $(u_1 \rho) (i, c, \mu)$ and $u_2 \rho$ are equal in all three
components and certainly $u_1, u_2 \in L \cap V$, from which it follows that
$(u_1, u_2) \in L_{(i, c, \mu)}$ as required.

Conversely, suppose $(u_1, u_2) \in L_{(i, c, \mu)} \cap (V \times V)$, and that
$u_1 \rho \pi_3 = \lambda$.  Then
$u_2 \rho = (u_1 \rho) (i, c, \mu)$, and using \eqref{eqn:autoreesb} and
equating second components, we obtain
$$u_2 \phi \sigma = u_2 \rho \pi_2
= (u_1 \rho \pi_2) P_{u_1 \rho \pi_3, i} c
= (u_1 \phi \sigma) P_{\lambda i} c
= (u_1 \phi \sigma) (b \sigma).$$
But now by the definition of $K_b$, it follows that $(u_1 \phi, u_2 \phi) \in K_b$.
It follows, as discussed above, that $K_b$ is synchronously regular, as
required to complete the proof.
\end{proof}

This result extends to cover the case where $S$ is only weakly
right-ideal-generated by a row cross-section of $P$.

\begin{theorem}\label{thm:auto_rees_mayimp_sgpoid2}
Let $M = M^0(S; F(I), G(\Lambda); P)$ be a Rees matrix semigroup with zero over
a semigroupoid $S$. Suppose $M$ is automatic, and the set of non-identity elements
of $S$ is weakly right-ideal-generated
by a row cross-section of $P$. Then $S$ is automatic.
\end{theorem}
\begin{proof}
Since $M$ is automatic it is certainly finitely generated, so
Theorem~\ref{thm:fgp_rees} tells
us that $S$ is finitely generated and that the indexing sets $I$ and
$\Lambda$ are finite. Let $\overline{S}$ be the category obtained from
$S$ by adjoining a new identity arrow at every object $e \in S^0$ which
does not already have one. Define $\overline{I}$,
$\overline{\Lambda}$, $\overline{P}$, $\overline{M}$ and
$\overline{U}$ just as in the proof of Theorem~\ref{thm:auto_sgpoid_imp_rees},
again noting that $\overline{I}$ and $\overline{\Lambda}$ are finite.
Still reasoning as in the proof of Theorem~\ref{thm:auto_sgpoid_imp_rees}, we
deduce that $M$ is a cofinite subsemigroupoid of $\overline{M}$. It
follows by Theorem~\ref{thm:auto_cofinite_subsgpoid} (or \cite[Theorem~1.1]{Hoffmann02})
that $\overline{M}$
is automatic.
We also deduce that $\overline{S} \setminus \overline{U}$ is finite.

Next, we wish to show that the set of non-identity elements of
$\overline{S}$ is strongly right-ideal-generated by a row cross-section of $\overline{P}$.
Since the non-identity elements of $S$ are weakly right-ideal-generated
by a row cross-section of $P$, we can choose some set $\Lambda' \subseteq \Lambda$ such
that the restriction of $G$ to $\Lambda'$ is bijective, and every non-identity
element of $S$ can be written in the form $P_{\lambda i} c$ or $P_{\lambda i}$
for some $\lambda \in \Lambda'$, $i \in I$ and (where appropriate) $c \in S$.
Let
$$\overline{\Lambda}' = \Lambda' \cup (\overline{\Lambda} \setminus \Lambda).$$

Now suppose $s$ is a non-identity element of $\overline{S}$. Notice that,
because we have adjoined an identity only where $S$ did not already have
one, every identity in $S$ remains an identity in $\overline{S}$. Thus,
we can assume that $s$ is a non-identity element of $S$. Thus, there exist $\lambda \in \Lambda'$ and $i \in I$ such that either
$s = P_{\lambda i} c$ for some $c \in S \subseteq \overline{S}$, or
$s = P_{\lambda i} = P_{\lambda i} e$ where $e \in \overline{S}$ is the identity at
$P_{\lambda i} \omega = i F$. Thus, the set of non-identity
elements of $\overline{S}$ is strongly right-ideal-generated by a
row cross-section of $P$.

It now follows by Theorem~\ref{thm:auto_rees_mayimp_sgpoid1} that
$\overline{S}$ is automatic.
But as in the proof of Theorem~\ref{thm:auto_sgpoid_imp_rees}, we
deduce that
$S$ is a cofinite subsemigroupoid of $\overline{S}$. Now by
Theorem~\ref{thm:auto_cofinite_subsgpoid},
it follows that $S$ is automatic, as required.
\end{proof}

\section{Rees Matrices and Prefix-Automaticity}\label{sec:pautorees}

In this section, we turn our attention to prefix-automaticity.
The following result, which provides a sufficient condition for a Rees
matrix semigroup over a semigroupoid to be prefix-automatic,
generalises \cite[Theorem~7.2]{Silva00b} in the case of prefix-automaticity.
As with Theorems~\ref{thm:auto_rees_mayimp_sgpoid1} and \ref{thm:auto_rees_mayimp_sgpoid2}, we shall prove this result
first with the hypothesis that the semigroupoid $S$ is a small category and that the non-identity
elements of $S$ are strongly right-ideal-generated by a row cross-section of the
sandwich matrix $P$. We
shall then extend the result to the more general case in which $S$ is a
semigroupoid with non-identity elements weakly right-ideal-generated by a
row cross-section of $P$.

\begin{theorem}\label{thm:pauto_cat_mayimp_rees}
Let $M = M^0(S; F(I), G(\Lambda); P)$ be a finitely generated Rees matrix semigroup with zero over
a small category $S$. If $S$ is prefix-automatic, and the non-identity elements
of $S$ are strongly right-ideal-generated by a row cross-section of $P$, then $M$ is prefix-automatic.
\end{theorem}
\begin{proof}
Let $\Lambda' \subseteq \Lambda$ be such that the restriction of $G$
to $\Lambda'$ is bijective, and every non-identity arrow in $S$ can be
written in the
form $P_{\lambda i} s$ for some $\lambda \in \Lambda'$, $i \in I$ and $s \in S$.

We define $U$ as in the proof of Theorem~\ref{thm:auto_cat_imp_rees}, once again
deducing from the fact that $M$ is finitely generated, that $I$,
$\Lambda$ and $S \setminus U$ are all finite and that $U$ is a subsemigroupoid
of $S$. It follows by Theorem~\ref{thm:auto_cofinite_subsgpoid} that $U$
is prefix-automatic, and by \cite[Proposition~4.4]{KambitesAutoCat}
 that $U$ has
a prefix-closed automatic structure. Note that since $S$ is finitely
generated, it has finitely many objects and hence finitely many local
identities. It follows by Proposition~\ref{prop:no_duplicate_letters}
that $U$ has a prefix-closed automatic structure $(X, K, \sigma)$ with the property
that the restriction of $\rho$ to $X$ is injective. (Note that, since we require our automatic
structure to be prefix-closed, we cannot insist that it should
also be a cross-section.)

We define $A$, $B$, $z$, $L$, $H$, $\rho$ and each $X_{i \lambda}$ and $\phi_{i \lambda}$ exactly
as in the proof of Theorem~\ref{thm:auto_cat_imp_rees}, and claim now that
$(A \cup B \cup \lbrace z \rbrace, L, \rho)$ is a prefix-automatic
structure for $S$.

By exactly the same arguments as in the proof of Theorem~\ref{thm:auto_cat_imp_rees},
we deduce that $(A \cup B \cup \lbrace z \rbrace, L, \rho)$ is a regular
choice of representatives for $M$, and that $L_a$ is synchronously regular for all
letters $a \in A \cup B \cup \lbrace z \rbrace$. However, since $(X, K, \sigma)$
is not here assumed to be a cross-section for $S$, we cannot deduce that
$(A \cup B \cup \lbrace z \rbrace, L, \rho)$ is a cross-section for $M$, and
we must work a little harder to show that $L_=$ is regular.

By Lemma~\ref{lem:limitedwork1}, it will suffice to show that
$L_= \cap (A^+ \times A^+)$
is synchronously regular.
Now certainly
$$L_= \cap (A^+ \times A^+) = \bigcup_{i \in I, \lambda \in \Lambda} L_{i, \lambda}$$
where
$$L_{i,\lambda} = \lbrace (u, v) \in L_= \cap (A^+ \times A^+) \mid u \rho \pi_1 = v \rho \pi_1 = i, u \rho \pi_3 = v \rho \pi_3 = \lambda \rbrace.$$
But
\begin{align*}
L_{i, \lambda} &= \lbrace (u, v) \in (K \cap X_{i \lambda}) \phi_{i \lambda} \times (K \cap X_{i \lambda}) \phi_{i \lambda} \mid u \rho = v \rho \rbrace \\
                  &= \lbrace (x \phi_{i \lambda}, y \phi_{i \lambda}) \mid (x, y) \in K_= \cap (X_{i \lambda} \times X_{i \lambda}) \rbrace.
\end{align*}
From our argument in the proof of Theorem~\ref{thm:auto_cat_imp_rees},
we know that $\phi_{i \lambda}$ is strongly synchronous regularity preserving,
 and we know
that $K_=$ is synchronously regular, so it follows that each $L_{i,\lambda}$
is synchronously regular. Hence, $L_= \cap (A^+ \times A^+)$ is a finite union of synchronously
regular binary relations, and by Lemma~\ref{lem:syncrat_prop}(iv) is synchronously
regular as required.

Thus, we conclude that $(A \cup B \cup \lbrace z \rbrace, L, \rho)$
is an automatic structure for $M$.
It remains only to show that the language
$$L'_= = \lbrace (u,v) \mid u \in L, v \in \OpPref(L), u \rho = v \rho \rbrace$$
is synchronously regular. First, observe that we can write
\begin{align*}
L'_= = &(L'_= \cap (B \times \OpPref(L))) \cup (L'_= \cap (L \times B)) \cup (L'_= \cap (\lbrace z \rbrace \times \OpPref(L)) \\
       &\cup (L'_= \cap (L \times \lbrace z \rbrace)) \cup (L'_= \cap (A^+ \times A^+)).
\end{align*}

Clearly, $z$ and letters in $B$ are unique representatives
in $\OpPref(L)$ for the respective elements they represent. Thus, we have
$$L'_= \cap (B \times \OpPref(L))) = L'_= \cap (L \times B) = \lbrace (b, b) \mid b \in B \rbrace$$
and
$$L'_= \cap (\lbrace z \rbrace \times \OpPref(L)) = L'_= \cap (L \times \lbrace z \rbrace) = \lbrace (z,z) \rbrace$$
so that the first four components of the union are synchronously regular.

It remains to show that
$L'_= \cap (A^+ \times A^+)$ is synchronously regular.
Define
$$\overline{L} = L \cup \lbrace a_{i 1_s 1_s \lambda} \mid i \in I, s \in S^0, \lambda \in \Lambda, i F = s = \lambda G \rbrace.$$
Clearly, since $I$ and $\Lambda$ are finite and $S$ has finitely many objects, $\overline{L} \setminus L$ is
finite. It follows by Lemma~\ref{lem:limitedwork2} that $\lbrace A \cup B \cup \lbrace z \rbrace, \overline{L}, \rho)$
is an automatic structure for $M$.

Given any letter $a = a_{i s t \lambda} \in A$, we define
$\overline{a} = a_{i s 1_{s \omega} \lambda'}$, where
$\lambda' = s \omega G'$ is the unique element in $\Lambda'$ satisfying
$\lambda' G = s \omega$.

Now consider a non-empty (not necessarily proper) prefix $a_1 \dots a_k$ of a word $a_1 \dots a_n$ in $L \cap A^+$.
We claim that $a_1 \dots a_{k-1} \overline{a_k} \in \overline{L}$.
If $k=1$ then $\overline{a_k}$ is of the form $a_{i 1_s 1_s \lambda}$, and
so is by definition in $\overline{L}$.
Otherwise, we observe that from the definitions of $L$ and $\phi_{i \lambda}$ we
must have
$$a_1 \dots a_n = (y_1 \dots y_{n-1}) \phi_{i \lambda}$$
for some $i \in I$, $\lambda \in \Lambda$ and some word $y_1 \dots y_{n-1} \in K \cap X_{i \lambda}$. But $K$ is
prefix-closed, so also $y_1 \dots y_{k-1} \in K$.
Furthermore, if we let $\mu \in \overline{\Lambda}$ with $\mu G = y_{k-1} \omega$
then we have $y_1 \dots y_{k-1} \in K \cap X_{i \mu}$ and
$$(y_1 \dots y_{k-1}) \phi_{i \mu} = a_1 \dots a_{k-1} \overline{a_k}.$$
Thus, $a_1 \dots a_{k-1} \overline{a_k} \in L \subseteq \overline{L}$ as claimed.

Recall from page \pageref{eqn:hdefine} the definitions of the subset
$H \subseteq S$ and the alphabet $A$. Let $s \in H$ and $\lambda \in \Lambda$ be such that $s \omega = \lambda F$.
We consider separately the cases
in which $s$ is and is not an identity element in $S$.

First, suppose that $s$ is not an identity element.
Then by assumption, we can choose $\lambda' \in \Lambda'$,
$i_s \in I$ and $t_s \in S$ such that $s = P_{\lambda' i_s} t_s$.
Then certainly $i_s F = t_s \alpha$ and
$t_s \omega = s \omega = \lambda F$, so there exists an element
$(i_s, t_s, \lambda) \in M$.

Since we have shown that $(A \cup B \cup \lbrace z \rbrace, \overline{L}, \rho)$ is
an automatic structure for $M$, we deduce using
\cite[Proposition~4.2]{KambitesAutoCat}
that the language
$$\overline{L}_{(i_s, t_s, \lambda)} = \lbrace (u, v) \in \overline{L} \times \overline{L} \mid (u \rho) (i_s, t_s, \lambda) = v \rho \rbrace$$
is synchronously regular. We define a new language
\begin{align*}
N_{s, \lambda} = \lbrace (u, a_1 \dots a_k) \mid &u \in L, a_1 \dots a_k \in \OpPref(L), a_k
 = a_{i, s', s, \lambda}, \\
&(a_1 \dots a_{k-1} \overline{a_k}, u) \in \overline{L}_{(i_s, t_s, \lambda)} \rbrace
\end{align*}
It follows easily from Proposition~\ref{prop:syncrat_changelastletter} and
Lemma~\ref{lem:syncrat_prop} that $N_{s,\lambda}$ is synchronously regular.

Now for any $u \in L$ and $a_1 \dots a_k \in \OpPref(L)$ with $a_k$ of the
form
$a_{i,s',s,\lambda}$ we have
$$\overline{a_k} \rho \pi_3 = s' \omega G' = s \alpha G' = P_{\lambda' i_s} \alpha G' = \lambda'$$
so that
\begin{align*}
(u, a_1 \dots a_k) \in N_{s,\lambda} &\iff (a_1 \dots a_{k-1} \overline{a_k}, u) \in \overline{L}_{(i_s, t_s, \lambda)} \\
&\iff (a_1 \dots a_{k-1} \overline{a_k}) \rho (i_s, t_s, \lambda) = u \rho \\
&\iff (a_1 \dots a_{k-1}) \rho (i, s' P_{\lambda' i_s} t_s, \lambda) = u \rho \\
&\iff (a_1 \dots a_{k-1}) \rho (i, s' s, \lambda) = u \rho \\
&\iff (a_1 \dots a_k) \rho = u \rho \\
&\iff (u, a_1 \dots a_k) \in L_='.
\end{align*}

Next, we consider the case in which $s$ is a local identity in $S$. In
this case, we define
\begin{align*}
N_{s, \lambda} = \lbrace (&b_1 \dots b_n, a_1 \dots a_k) \mid
 a_1 \dots a_k \in \OpPref(L), b_1 \dots b_n \in L, \\
&a_k = a_{i, s', s, \lambda},
b_n = a_{j, t, 1_{t \alpha}, \lambda},
 (a_1 \dots a_{k-1} \overline{a_k}, b_1 \dots b_{n-1} \overline{b_n})
 \in \overline{L}_= \rbrace.
\end{align*}

Since $(A \cup B \cup \lbrace z \rbrace, \overline{L}, \rho)$ is an automatic
structure for $M$, we deduce that $\overline{L}_=$ is regular, and then by
Proposition~\ref{prop:syncrat_changelastletter} and Lemma~\ref{lem:syncrat_prop}(i)
that $N_{s,\lambda}$ is synchronously regular.

Now for any $a = a_1 \dots a_k \in \OpPref(L)$ and $b = b_1 \dots b_n \in L$ with
$a_k$ of the form $a_{i,s',s,\lambda}$ and $b_n$ of the form
$a_{j, t, 1_{t \alpha}, \lambda}$ we have
\begin{align*}
(b, a) \in N_{s,\lambda} &\iff
 (a_1 \dots a_{k-1} \overline{a_k}, b_1 \dots b_{n-1} \overline{b_n}) \in \overline{L}_= \\
&\iff (a_1 \dots a_{k-1} \overline{a_k}) \rho = (b_1 \dots b_{n-1} \overline{b_n}) \rho \\
&\iff (a_1 \dots a_{k-1}) \rho (i, s', \lambda G G') = (b_1 \dots b_{n-1}) \rho (j, t, \lambda G G') \\
&\iff (a_1 \dots a_{k-1}) \rho (i, s's, \lambda) = (b_1 \dots b_{n-1}) \rho (j, t 1_{t \alpha}, \lambda) \\
&\iff (a_1 \dots a_k) \rho = (b_1 \dots b_n) \rho \\
&\iff (b_1 \dots b_n, a_1 \dots a_k) \in L_='.
\end{align*}

It follows from the two cases considered that
$$L'_= = \bigcup_{s \in H, \lambda \in \Lambda} N_{s,\lambda}$$
so that $L'_=$ is a union of finitely many synchronously regular
languages, and hence by Lemma~\ref{lem:syncrat_prop}(iv) is
synchronously regular.
\end{proof}

We now extend this result as described above.

\begin{theorem}\label{thm:pauto_sgpoid_mayimp_rees}
Let $M = M^0(S; F(I), G(\Lambda); P)$ be a finitely generated Rees matrix semigroup with zero over
a semigroupoid $S$. If $S$ is prefix-automatic, and the set of non-identity
elements in $S$ is weakly right-ideal-generated by a row cross-section of $P$, then
$M$ is prefix-automatic.
\end{theorem}

\begin{proof}
We combine the methods used to prove Theorems~\ref{thm:auto_sgpoid_imp_rees} and
\ref{thm:auto_rees_mayimp_sgpoid2}.

As in the proof of Theorem~\ref{thm:auto_sgpoid_imp_rees}, we deduce that
the indexing sets $I$ and $\Lambda$ are finite. We define $\overline{S}$, $\overline{I}$,
$\overline{\Lambda}$, $\overline{P}$, $\overline{M}$ and
$\overline{U}$ as in the proof of Theorem~\ref{thm:auto_rees_mayimp_sgpoid2},
noting that $\overline{S}$ has an adjoined identity only where there was
not already an identity in $S$.

As in Theorem~\ref{thm:auto_rees_mayimp_sgpoid2}, we
deduce that $S$ is a cofinite subsemigroupoid of $\overline{S}$, that $M$
is a cofinite subsemigroup of $\overline{M}$, and that the non-identity
elements of $\overline{S}$ are strongly right-ideal-generated by a row
cross-section of $\overline{P}$.

It follows by Theorem~\ref{thm:auto_cofinite_subsgpoid}
that $\overline{S}$ is prefix-automatic, by
Theorem~\ref{thm:pauto_cat_mayimp_rees} that $\overline{M}$ is
prefix-automatic and then by Theorem~\ref{thm:auto_cofinite_subsgpoid}
that $M$ is automatic, as required.
\end{proof}

We now show that prefix-automaticity in a Rees matrix
semigroup is a sufficient condition for prefix-automaticity in the underlying
semigroupoid. The following result generalises \cite[Theorem~4.2]{Descalco01}.

\begin{theorem}\label{thm:pauto_rees_imp_sgpoid}
Let $M = M^0(S; F(I), G(\Lambda); P)$ be a Rees matrix semigroup with zero over
a semigroupoid $S$. If $M$ is prefix-automatic then $S$ is
prefix-automatic.
\end{theorem}
\begin{proof}
Suppose $M$ is prefix-automatic. Then by \cite[Proposition~4.4]{KambitesAutoCat}
$M$ has a prefix-closed automatic structure, and by Proposition~\ref{prop:no_duplicate_letters},
$M$ has a prefix-closed automatic structure $(A, L, \rho)$ with the property that the
restriction of $\rho$ to $A$ is injective. (Once again, we note that, because
we require an automatic structure which is prefix-closed, we cannot ask also
that it be a cross-section.)

Choose some subsets $I' \subseteq I$ and $\Lambda' \subseteq \Lambda$ such
that $F$ and $G$ restrict to bijections on $I'$ and $\Lambda'$ respectively.
Define $F'$, $G'$, $C$, $D$, $X$, $\sigma$, $V$, $\phi$ and $K$ exactly as in
the proof
of Theorem~\ref{thm:auto_rees_mayimp_sgpoid1}.
Reasoning as before, we deduce
that $\phi$ is strongly regularity preserving and strongly synchronous regularity preserving,
and that for any $w \in V$ we have
\begin{equation}\label{eqn:autoreesb_repeat}
w \rho = (w \phi \sigma \alpha F', w \phi \sigma, w \phi \sigma \omega G').
\end{equation}
We deduce also that $(X, K, \sigma)$ is a regular choice of representatives for $S$, although
we can no longer conclude that it is a cross-section.
We claim that $(X, K, \sigma)$ is a prefix-automatic structure for $S$.

Because of the limited role played by the alphabet in the definition of
an automatic structure, we can assume without loss of generality that
for every letter $a \in A$, there is a letter in $A$ representing the
element
$$(a \rho \pi_1, a \rho \pi_2, a \rho \pi_3 G G') \in M.$$ For each
$a \in A$, let $\overline{a} \in A$ be such a letter.

Let $c \in X^1 \cap X^0$ be an edge or a path of length 0 in $X$. Consider
the binary relation
$$K_c' = \lbrace (x, y) \in \OpPref(K) \times K \mid x \omega = c \alpha, (xc) \sigma = y \sigma \rbrace.$$
Our aim is to show that $K_c'$ is synchronously regular. We claim that
\begin{equation}\label{eqn:kcclaim}
K_c' = \left( \bigcup_{\lambda \in \Lambda, i \in I, b \in C \cup CD} K_{\lambda,i,b} \right) \cup \left( \bigcup_{x \in \OpPref(K) \cap X^{\leq 4}} K_c' \cap ( \lbrace x \rbrace \times K) \right)
\end{equation}
where each
\begin{align*}
K_{\lambda,i,b} = \lbrace & (((w_1 \dots w_{n-1} \overline{w_n}) \phi) d_{\lambda i} b, v \phi) \mid w_1, w_2, \dots, w_n \in A, v \in L \cap V \\
&w = w_1 \dots w_n \in L, w \rho \pi_3 = \lambda, (w_1 \dots w_n, v) \in L_{(i,(bc) \sigma, c \omega G')} \rbrace.
\end{align*}
To prove the claim, suppose first that $(x,y) \in K_c'$. Clearly, if $|x| \leq 4$
then, since $(x,y) \in K_c' \cap ( \lbrace x \rbrace \times K)$,
we see that $(x,y)$ is contained in the right-hand-side of \eqref{eqn:kcclaim}
and we are done.

Now
suppose $|x| \geq 5$. Certainly $y \in K = (L \cap V) \phi$, so we can
write $y = v \phi$ for some $v \in L \cap V$. Also $x \in \OpPref(K) = \OpPref((L \cap V) \phi)$,
so certainly we can choose $w \in L \cap V$ and $z \in X^*$ with $xz = w \phi$.
From the definition of $\phi$, we have $xz \in C(DC)^*$. Since
$|x| \geq 5$, it follows
that we can write $x = a d_{\lambda i} b$ for some $a \in CDC(DC)^*$,
$\lambda \in \Lambda$, $i \in I$ and either $b \in C$ or $b \in CD$.

Suppose $w = w_1 \dots w_k$.  Then it is easily verified from the definition
of $\phi$ (see page \pageref{eqn:phidefine}) that
$a = (w_1 \dots w_{n-1} \overline{w_n}) \phi$
where $n = (|a|+1) / 2$.
Let $u = w_1 \dots w_n$, noticing that we have $u \rho \pi_3 = \lambda$.
Now we have
\begin{align*}
(u \rho) (i, (bc) \sigma, c \omega G')
&= (u \rho \pi_1, u \rho \pi_2, \lambda) (i, (bc) \sigma, c \omega G') \\
&= (u \rho \pi_1, a \sigma, \lambda) (i, (bc) \sigma, c \omega G') \\
&= (u \rho \pi_1, (a \sigma) P_{\lambda i} (bc) \sigma, c \omega G') \\
&= (u \rho \pi_1, (a d_{\lambda i} b c) \sigma, c \omega G') \\
&= (v \rho \pi_1, (x c) \sigma, v \rho \pi_3) \\
&= (v \rho \pi_1, y \sigma, v \rho \pi_3) \\
&= v \rho
\end{align*}
so that $(u, v) \in L_{(i, (bc) \sigma, c \omega G')}$
and hence $\left( ((w_1 \dots w_{n-1} \overline{w_n}) \phi \right) d_{\lambda i} b, v \phi) \in K_{\lambda,i,b}$.

Conversely, suppose $w = w_1 \dots w_n \in L$ with
$$\left( (w_1 \dots w_{n-1} \overline{w_n}) \phi \right) d_{\lambda i} b, v \phi) \in K_{\lambda,i,b}$$
for some $\lambda \in \Lambda, i \in I$ and $b \in C \cup CD$. Let
$a = (w_1 \dots w_{n-1} \overline{w_n}) \phi.$
Now
\begin{align*}
(v \rho \pi_1, (((w_1 \dots w_{n-1} \overline{w_n}) \phi) d_{\lambda i} b c) \sigma, v \rho \pi_3)
&= (w \rho \pi_1, (a d_{\lambda i} b c) \sigma, c \omega G') \\
&= (w \rho \pi_1, (a \sigma) P_{\lambda i} (bc) \sigma, c \omega G') \\
&= (w \rho \pi_1, a \sigma, \lambda) (i, (bc) \sigma, c \omega G') \\
&= (w \rho \pi_1, w \rho \pi_2, \lambda) (i, (bc) \sigma, c \omega G') \\
&= (w \rho) (i, (bc \sigma), c \omega G') \\
&= v \rho
\end{align*}
so in particular we have
$$(((w_1 \dots w_{n-1} \overline{w_n}) \phi) d_{\lambda i} b \sigma) (c \sigma) = v \phi \sigma$$
so that
$(((w_1 \dots w_{n-1} \overline{w_n}) \phi) d_{\lambda i} b \sigma, v \phi \sigma) \in K'_c$
as required.

Also, it is clear that if $x \in \OpPref(K) \cap X^{\leq 4}$ and $y \in K$ are
such that $(x,y) \in K_c' \cap ( \lbrace x \rbrace \times K)$
then we must have $(x,y) \in K_c'$. Thus, we have justified our claim
that \eqref{eqn:kcclaim} holds.

Now using \cite[Proposition~4.2]{KambitesAutoCat}
we see that each binary relation of the form
$L_{(i, (bc) \sigma, c \omega G')}$ is synchronously regular.
Furthermore, the functions given by
$$A \to A, a \mapsto \overline{a}, \text{ and}$$
$$X \to X^{|b| + 2}, x \mapsto x d_{\lambda i} b$$
both satisfy the conditions of Proposition~\ref{prop:syncrat_changelastletter}.
Since $\phi$ is
also synchronous regularity preserving, it follows that each
$K_{\lambda,i,b}$ is synchronously regular.

Now let $x \in \OpPref(K) \cap X^{\leq 4}$. Choose a word $u$ in $L$
representing the element $(x \alpha F', (x \sigma) (c \sigma), c \omega G')$. Now
for any word $v \phi \in (L \cap V) \phi = K$, we have $(x c) \sigma = v \phi \sigma$
if and only if $u \rho = v \rho$, which in turn is true exactly if
$(u, v) \in L_=$. Thus, we have
$$K_c' \cap ( \lbrace x \rbrace \times K) = \lbrace x \rbrace \times ((\lbrace u \rbrace \times L) \cap L_=) \pi_2 \phi.$$
We know that $L$ is regular and that $L_=$ is synchronously regular, so
it follows by Lemma~\ref{lem:syncrat_prop}(ii), (iii) and (iv)
and the fact that $\phi$ is strongly synchronous regularity preserving
that $K_c' \cap (\lbrace x \rbrace \times K)$ is synchronously regular.

We have shown that each $K_c'$ is a finite union of synchronously regular
binary relations, and it follows by Lemma~\ref{lem:syncrat_prop}(iv) that each
$K_c'$ is synchronously regular.
Now for every $c \in X^1 \cup X^0$ we see that
$K_c = K_c' \cap (K \times K)$
is synchronously regular by Lemma~\ref{lem:syncrat_prop}(iv), so that
$(X, L, \sigma)$ is an automatic structure for $S$. Moreover,
$$K_=' = \bigcup_{c \in X^0} (K_c')^{-1}$$
is also synchronously regular by Lemma~\ref{lem:syncrat_prop}(iv), so
that $(X, K, \sigma)$ is a prefix-automatic structure for $S$, as required.
\end{proof}

\section{Rees Matrix Semigroups Without Zero}\label{sec:reeswithoutzero}

Theorem~\ref{thm:auto_cofinite_subsgpoid} ensures that our results about Rees matrix
semigroups with zero adapt easily to the case of Rees matrix semigroups without
zero.

\begin{theorem}\label{thm:auto_sgpoid_imp_reeswithoutzero}
Let $M = M(S; F(I), G(\Lambda); P)$ be a Rees matrix semigroup (without zero) over a
semigroupoid $S$. Suppose $S$ is automatic and $M$ is finitely
generated. Then $M$ is automatic.
\end{theorem}
\begin{proof}
Consider the Rees matrix semigroup with zero
$$M' = M^0(S; F(I), G(\Lambda); P).$$
Then $M = M' \setminus \lbrace 0 \rbrace$ is a cofinite subsemigroup
of $M'$. So by \cite[Theorem~1.1]{Ruskuc98}, $M'$ is finitely generated. Now by
Theorem~\ref{thm:auto_sgpoid_imp_rees}, $M'$ is automatic, and by
Theorem~\ref{thm:auto_cofinite_subsgpoid} (or \cite[Theorem~1.1]{Hoffmann02})
it follows that $M$ is automatic.
\end{proof}

\begin{theorem}\label{thm:auto_reeswithoutzero_mayimp_sgpoid}
Let $M = M(S; F(I), G(\Lambda); P)$ be a Rees matrix semigroup (without zero) over
a semigroupoid $S$. Suppose the set of non-identity elements of $S$ is
weakly right-ideal-generated by a
row cross-section of $P$. If $M$ is automatic then $S$ is automatic.
\end{theorem}
\begin{proof}
Consider the Rees matrix semigroup with zero
$$M' = M^0(S; F(I), G(\Lambda); P).$$
Then $M = M' \setminus \lbrace 0 \rbrace$ is a cofinite subsemigroupoid of $M'$.
Now if $M$ is automatic then by Theorem~\ref{thm:auto_cofinite_subsgpoid},
$M'$ is automatic, and by Theorem~\ref{thm:auto_rees_mayimp_sgpoid2},
$S$ is automatic.
\end{proof}

\begin{theorem}\label{thm:pauto_sgpoid_mayimp_reeswithoutzero}
Let $M = M(S; F(I), G(\Lambda); P)$ be a finitely generated Rees matrix
semigroup (without zero) over a semigroupoid $S$. If $S$
is prefix-automatic, and the set of non-identity elements of $S$ is weakly
right-ideal-generated by a row cross-section of $P$,
then $M$ is prefix-automatic.
\end{theorem}
\begin{proof}
Consider the Rees matrix semigroup with zero
$$M' = M^0(S; F(I), G(\Lambda); P).$$
Then $M = M' \setminus \lbrace 0 \rbrace$ is a cofinite subsemigroup of $M'$. So
by \cite[Theorem~1.1]{Ruskuc98}, $M'$ is finitely generated. Now by
Theorem~\ref{thm:pauto_sgpoid_mayimp_rees}, $M'$ is prefix-automatic, and by
Theorem~\ref{thm:auto_cofinite_subsgpoid} it follows that $M$ is prefix-automatic.
\end{proof}

\begin{theorem}\label{thm:pauto_reeswithoutzero_imp_sgpoid}
Let $M = M(S; F(I), G(\Lambda); P)$ be a Rees matrix semigroup (without zero) over
a semigroupoid $S$. If $M$ is prefix-automatic then $S$ is
prefix-automatic.
\end{theorem}
\begin{proof}
Consider the Rees matrix semigroup with zero
$$M' = M^0(S; F(I), G(\Lambda); P).$$
Then $M = M' \setminus \lbrace 0 \rbrace$ is a cofinite subsemigroupoid of $M'$.
Now if $M$ is prefix-automatic then by Theorem~\ref{thm:auto_cofinite_subsgpoid}, $M'$
is prefix-automatic. By Theorem~\ref{thm:pauto_rees_imp_sgpoid}, it follows
that $S$ is prefix-automatic.
\end{proof}

\section{Closing Remarks}\label{sec:remarks}

The curious relationship between the results of Section~\ref{sec:pautorees}
and those of Section~\ref{sec:autorees} seems to demand comment. In the case
of automaticity, showing that automaticity in a Rees matrix semigroup is a
\textit{sufficient} condition for automaticity in the underlying semigroupoid
requires a right-ideal-generation condition, while the converse implication
does not. In the case of prefix-automaticity, the situation is entirely
reversed -- the right-ideal-generation condition is required only to show
that prefix-automaticity in a Rees matrix semigroup is a \textit{necessary}
condition for prefix-automaticity in the underlying semigroupoid.

In \cite{KambitesAutoCat} we remarked that it is an open question whether
every automatic semigroup is prefix-automatic.
Theorem~\ref{thm:cons_iff_sgpoid} implies that this question
is no harder in the ostensibly more general semigroupoid context, that is, that every automatic
semigroupoid is prefix-automatic, exactly if every automatic semigroup is prefix-automatic.
However, the following questions do naturally arise.
\begin{question}\label{qn:existautorees_over_nonauto_sgp}
Does there exist an automatic Rees matrix semigroup over a non-automatic semigroup?
\end{question}
\begin{question}\label{qn:existautorees_over_nonauto_sgpoid}
Does there exist an automatic Rees matrix semigroup over a non-automatic semigroupoid?
\end{question}
Clearly, if every automatic semigroup (and hence every automatic semigroupoid) is
prefix-automatic, then Theorem~\ref{thm:pauto_rees_imp_sgpoid} gives a negative
answer to both of these questions. However, if there \textit{are} automatic semigroups
(and hence semigroupoids) which are not prefix-automatic, then the answers to
Questions~\ref{qn:existautorees_over_nonauto_sgp} and
\ref{qn:existautorees_over_nonauto_sgpoid} could be both positive, both negative, or
negative and positive respectively.

\section*{Acknowledgements}
This paper was written while the author was at Carleton University,
supported by the Leverhulme Trust. The research documented was conducted
while the author was a research student at the University of York, funded
by an EPSRC Doctoral Studentship. The author would like to thank John
Fountain for all his advice and guidance, as well as Vicky Gould, Nik
Ru\v{s}kuc and the anonymous referees for many helpful comments. He
would also like to thank Kirsty for all her support and encouragement.

\bibliographystyle{plain}

\end{document}